\documentclass[reqno,12pt]{amsart}
\usepackage{amsmath,amsthm,amssymb,amsfonts,amscd,graphicx}
\input xy
\xyoption{all}

\theoremstyle{plain} 
	\newtheorem{thm}{Theorem}[section]
	\newtheorem*{thm*}{Theorem}
	
	\newtheorem{cora}[thm]{${\rm \bf Corollary}^{*}$}
	
        \newtheorem{lem}[thm]{Lemma}
	\newtheorem{lema}[thm]{${\rm \bf Lemma}^{*}$}
	\newtheorem{lemd}[thm]{${\rm \bf Lemma}^{\dagger}$}
	
	\newtheorem{sublema}[thm]{${\rm \bf Sublemma}^{*}$}
	\newtheorem{sublemd}[thm]{${\rm \bf Sublemma}^{\dagger}$}
	\newtheorem{prop}[thm]{Proposition}
	\newtheorem*{prop*}{Proposition}
	\newtheorem{propa}[thm]{${\rm \bf Proposition}^{*}$}
	\newtheorem{propd}[thm]{${\rm \bf Proposition}^{\dagger}$}
	\newtheorem{conj}[thm]{Conjecture}
	\newtheorem*{conj*}{Conjecture}
	
\theoremstyle{definition}
	\newtheorem{defn}[thm]{Definition}

\theoremstyle{remark}
	\newtheorem{rem}[thm]{Remark}
	
	\newtheorem*{pf}{Proof}
\numberwithin{equation}{section}
\def\CC{{\mathbb C}}

\def\NN{{\mathbb N}}
\def\PP{{\mathbb P}}
\def\QQ{{\mathbb Q}}

\def\ZZ{{\mathbb Z}}

\def\F{{\mathcal F}}
\def\G{{\mathcal G}}

\def\I{{\mathcal I}}

\def\M{{\mathcal M}}

\def\O{{\mathcal O}}

\def\T{{\mathcal T}}
\def\X{{\mathcal X}}

\def\p{\partial }

\def\ns{{\nabla}\hspace{-1.4mm}\raisebox{0.3mm}{\text{\footnotesize{\bf /}}}}

\begin{document}
\title[On Frobenius Manifolds form Gromov--Witten Theory]
{On Frobenius Manifolds from Gromov--Witten Theory of Orbifold Projective Lines
with $r$ orbifold points}
\date{\today}
\author{Yuuki Shiraishi}
\address{Department of Mathematics, Graduate School of Science, Osaka University, 
Toyonaka Osaka, 560-0043, Japan}
\email{sm5021sy@ecs.osaka-u.ac.jp}
\begin{abstract}
We prove that the Frobenius structure constructed from the Gromov-Witten theory  
for an orbifold projective line with at most $r$ orbifold points is uniquely determined by the WDVV equations with certain 
natural initial conditions. 
\end{abstract}
\maketitle
\section{Introduction}

The (formal) Frobenius manifold is a certain complex (formal) manifold endowed with the Frobenius algebra structure on its tangent sheaf
whose product, unit, non--degenerate bilinear form and grading operator called the Euler vector field satisfy the special properties 
(for its definition and important properties, see Section~\ref{section2}).  
This structure was originally discovered by K. Saito in his study of primitive forms and their period mappings
on the deformation theory of isolated hypersurface singularities (\cite{st:1} and references therein) and was
rediscovered and formulated by Dubrovin \cite{d:1} in order to give coordinate--free expression for a solution of
the Witten--Dijkgraaf--Verlinde--Verlinde $($WDVV$)$ equations considered in two dimensional
topological field theories. Namely the Frobenius manifold can be also obtained from the Gromov--Witten theory for manifolds or orbifolds.
Here the Gromov--Witten theory for orbifolds by Abramovich--Graber--Vistoli \cite{agv:1} and Chen--Ruan \cite{cr:1} 
is summarized briefly as follows;  
Let $\X$ be an orbifold $($or a smooth proper Deligne--Mumford stack over $\CC)$.
Then, for non--negative integers $g, n\in \ZZ_{\ge 0}$ 
and $\beta\in H_2(X,\ZZ)$ where $X$ is the coarse moduli space of $\X$, 
the moduli space (stack) $\overline{\M}_{g,n}(\X,\beta)$ 
of orbifold (twisted) stable maps of genus $g$ 
with $n$-marked points of degree $\beta$ is defined.
There exists a virtual fundamental class $\left[\overline{\M}_{g,n}(\X,\beta)\right]^{\rm vir}$ 
and Gromov--Witten invariants of genus $g$ 
with $n$-marked points of degree $\beta$ are defined as usual except for 
that we have to use the orbifold cohomology group $H^*_{orb}(\X,\QQ)$:
\[
\left<\Delta_1,\dots, \Delta_n\right>_{g,n,\beta}^\X:=
\int_{\left[\overline{\M}_{g,n}(\X,\beta)\right]^{\rm vir}}{\rm ev}_1^*\Delta_1\wedge \dots \wedge 
{\rm ev}_n^*\Delta_n,\quad \Delta_1,\dots,\Delta_n\in H^*_{orb}(\X,\QQ),
\]
where we denote by ${\rm ev}^*_i:H^*_{orb}(\X,\QQ)\longrightarrow H^*(\overline{\M}_{g,n}(\X,\beta),\QQ)$
the induced homomorphism by the evaluation map.
We also consider the generating function (or formal power series)
\[
\F_g^\X:=\sum_{n,\beta}\frac{1}{n!}\left<{\bf t},\dots, {\bf t}\right>_{g,n,\beta}^\X q^\beta,
\quad {\bf t}=\sum_{i}t_i\Delta_i
\]
and call it the genus $g$ potential where $\{\Delta_i\}$ denotes a $\QQ$-basis of 
$H^*_{orb}(\X,\QQ)$.
The main result in \cite{agv:1, cr:1} tells us that the point axiom, the divisor axiom for a class in $H^2(\X,\QQ)$
and the associativity of the quantum product hold same as the Gromov--Witten theory for a usual manifold (see \cite{agv:1, cr:1} for details of these axioms).
In particular, the associativity of the quantum product implies the WDVV equations and it 
gives a formal Frobenius manifold $M$ whose structure sheaf $\O_M$, tangent sheaf $\T_M$ and Frobenius potential are defined as 
the algebra $\Lambda[[H^*_{orb}(\X,\CC)]]$ of formal power series in dual coordinates $\{t_i\}$ 
of the $\QQ$-basis $\{\Delta_i\}$ of $H^*_{orb}(\X,\QQ)$ over the Novikov field $\Lambda$ 
(roughly speaking, $\Lambda$ is the $\CC$-algebra of formal Laurent series in $q^{\beta_1},\dots, q^{\beta_\rho}$ where 
$\beta_1$, $\dots$, $\beta_\rho$ are effective $1$-cycles which generate the Mori cone of $X$),  
$\T_M:=H^*_{orb}(\X,\CC)\otimes_\CC\O_M$ and the genus zero 
potential $\F_0^\X$ respectively.
Let $r\ge 3$ be a positive integer.
Let $A$ be a multiplet $(a_1,a_2,\dots, a_r)$ of positive integers such that $2\le a_1\le a_2\le \dots \le a_r$ 
and $\Lambda=(\lambda_1,\ldots,\lambda_r)$ a multiplet of pairwise distinct elements of $\PP^{1}(\CC)$ normalized such that 
$\lambda_1=\infty$, $\lambda_2=0$ and $\lambda_3=1$.
Set $\mu_A=2+\sum_{k=1}^{r}(a_{k}-1)$ and $\chi_A:=2+\sum_{k=1}^{r}(-1+1/a_k)$.
We shall consider the orbifold projective line with $r$--orbifold points at $\lambda_1,\dots,\lambda_r$ whose orders are 
$a_1, a_2, \dots, a_r$, which is denoted by $\PP^{1}_{A,\Lambda}$ (see Definition~\ref{defn:gl}).
Here the number $\mu_A$ is regarded as the total dimension of the orbifold cohomology group $H^*_{orb}(\PP^{1}_{A,\Lambda},\CC)$ and the number  
$\chi_A$ is regarded as the orbifold Euler number of $\PP^{1}_{A,\Lambda}$. 
The main purpose of the present paper is to show that the Frobenius manifold $M^{GW}_{\PP^{1}_{A,\Lambda}}$ 
constructed from the Gromov--Witten theory for $\PP^{1}_{A,\Lambda}$ 
can be determined by the WDVV equations with certain natural initial conditions. Then we shall show the following uniqueness theorem 
which is our main result in the present paper and the natural generalization of the one in our previous paper \cite{ist:1}:
\begin{thm*}[Theorem~\ref{first}]
There exists a unique Frobenius manifold $M$ of rank $\mu_A$ and dimension one with flat coordinates 
$(t_1,t_{1,1},\dots ,t_{i,j},\dots ,t_{r,a_r-1},t_{\mu_A})$ satisfying the following conditions$:$
\begin{enumerate}
\item 
The unit vector field $e$ and the Euler vector field $E$ are given by
\[
e=\frac{\p}{\p t_1},\ E=t_1\frac{\p}{\p t_1}+\sum_{i=1}^{r}\sum_{j=1}^{a_i-1}\frac{a_i-j}{a_i}t_{i,j}\frac{\p}{\p t_{i,j}}
+\chi_A\frac{\p}{\p t_{\mu_A}}.
\]
\item 
The non--degenerate symmetric bilinear form $\eta$ on $\T_M$ satisfies
\begin{align*}
&\ \eta\left(\frac{\p}{\p t_1}, \frac{\p}{\p t_{\mu_A}}\right)=
\eta\left(\frac{\p}{\p t_{\mu_A}}, \frac{\p}{\p t_1}\right)=1,\\ 
&\ \eta\left(\frac{\p}{\p t_{i_1,j_1}}, \frac{\p}{\p t_{i_2,j_2}}\right)=
\begin{cases}
\frac{1}{a_{i_1}}\quad i_1=i_2\text{ and }j_2=a_{i_1}-j_1,\\
0 \quad \text{otherwise}.
\end{cases}
\end{align*}
\item 
The Frobenius potential $\F_{A}$ satisfies $E\F_{A}|_{t_{1}=0}=2\F_{A}|_{t_{1}=0}$,
\[
\left.\F_{A}\right|_{t_1=0}\in\CC\left[[t_{1,1}, \dots, t_{1,a_1-1}, 
\dots, t_{i,j},\dots, t_{r,1}, \dots, t_{r,a_r-1},e^{t_{\mu_A}}]\right].
\]
\item Assume the condition {\rm (iii)}. we have
\begin{equation*}
\F_{A}|_{t_1=e^{t_{\mu_A}}=0}=\sum_{i=1}^{r}\G_{A}^{(i)}, \quad \G_{A}^{(i)}\in \CC[t_{i,1},\dots, t_{i,a_i-1}],\ i=1,\dots,r.
\end{equation*}
\item 
Assume the condition {\rm (iii)}. In the frame $\frac{\p}{\p t_1}, \frac{\p}{\p t_{1,1}},\dots, 
\frac{\p}{\p t_{r,a_r-1}},\frac{\p}{\p t_{\mu_A}}$ of $\T_M$,
the product $\circ$ can be extended to the limit $t_1=t_{1,1}=\dots=t_{r,a_r-1}=e^{t_{\mu_A}}=0$.
The $\CC$-algebra obtained in this limit is isomorphic to
\[
\CC[x_1,x_2,\dots, x_r]\left/\left(x_ix_j, \ a_ix_i^{a_i}-a_jx_j^{a_j}
\right)_{1\le i\ne j\le r}\right.,
\]
where $\p/\p t_{i,j}$ are mapped to
$x^{j}_i$ for $i=1,\dots,r, j=1,\dots, a_{i}-1$ and $\p/\p t_{\mu_A}$ are mapped to $a_{1}x_{1}^{a_1}$.
\item The term 
\[
\displaystyle\left(\prod_{i=1}^{r}t_{i,1}\right)e^{t_{\mu_A}}
\]
occurs with the coefficient $1$ in $\F_{A}$. 
\end{enumerate}
\end{thm*}
Here we have two important results concerning the condition {\rm (iv)} in Theorem~\ref{first}.
First, the polynomial $\G_{A}^{(i)}$ in the condition {\rm (iv)} can be expressed by
the Frobenius potential $\F_{A_{i}}(t'_{1}, {\bf t'}_{3}, e^{t_{\mu_{A_{i}}}})$ of the Frobenius manifold $M_{A_{i}}$
in Theorem \ref{first} in \cite{ist:1} where $A_{i}=(1,1,a_{i})$ with $a_{i}\ge 2$
and $(t'_{1}, {\bf t'}_{3}, t_{\mu_{A_{i}}}):=(t'_{1},t'_{3,1},\dots,t'_{3,a_{i}-1}, t_{\mu_{A_{i}}})$ is the flat coordinate for $M_{A_{i}}$: 
\begin{prop*}[Proposition~\ref{atype decomp-prop}]
For the polynomial $\G_{A}^{(i)}$ in the condition {\rm (iv)} in Theorem \ref{first}, we have 
\begin{equation*}
\G_{A}^{(i)}=\F_{A_{i}}(0,{\bf t}_{i},0),
\end{equation*}
where ${\bf t}_{i}:=(t_{i,1},\dots,t_{i,a_{i}-1})$ is the $i$--th parts of the flat coordinate in Theorem~\ref{first}.
\end{prop*}
Second, we can derive the condition {\rm (iv)} from other conditions if the multiplet $A$ satisfies $2\le a_{2}<a_{3}$ (called ``general multiplet")
or under some weak condition if the multiplet $A$ satisfies $2=a_{1}=a_{2}<a_{3}$ (called ``semi--general multiplet").
This is a generalization of Proposition 3.24 in \cite{ist:1}: 
\begin{thm*}[Theorem \ref{sep}]
Suppose that $A$ is a general multiplet.
For a non--negative $\beta\in \ZZ^{\mu_A-2}$, we have
\[
c\left(\beta+\sum_{k=1}^3e_{i_k,j_{k}},0\right)\ne 0
\]
only if $i_1=i_2=i_3$.
Suppose that $A$ is a semi--general multiplet.
For a non--negative $\beta\in \ZZ^{\mu_A-2}$, we have
\[
c\left(\beta+\sum_{k=1}^3e_{i_k,j_{k}},0\right)\ne 0
\]
only if $i_1=i_2=i_3$ under the following condition:
\begin{itemize}
\item[{\rm (iv')}] If $a_{i_1}=a_{i_2}$ for some $i_1,i_2\in\{1,\dots,r\}$, then the
Frobenius potential $\F$ is invariant under
the permutation of parameters $t_{i_1,j}$ and $t_{i_2,j}$ $(j=1,\dots, a_{i_1}-1)$.
\end{itemize}
\end{thm*}
As a corollary of Theorem \ref{first}, the Frobenius structure constructed from the Gromov--Witten theory for $\PP^{1}_{A,\Lambda}$ can be uniquely 
reconstructed by the conditions in Theorem~\ref{first}:  
\begin{thm*}[Theorem~\ref{second}]
The conditions in Theorem~\ref{first} are satisfied by the Frobenius structure constructed
from the Gromov--Witten theory for $\PP^{1}_{A,\Lambda}$. 
\end{thm*}
In the rest of the present paper, we investigate the Frobenius potential which satisfies the same conditions with Theorem~\ref{first} 
except for the condition {\rm (vi)}. In other theory like the invariant theory of extended cuspidal Weyl groups 
which is expected as a mirror partner of $\PP^{1}_{A,\Lambda}$ (\cite{Sh-T:1}), the representation theoretic
meaning of this coefficient is not known at all and it is hard even to verify whether this coefficient is non--zero or not.
If this coefficient were non--zero, we can apply Theorem~\ref{first} to showing the isomorphism between the Frobenius manifold 
constructed from the invariant theory of an extended cuspidal Weyl group and $M^{GW}_{\PP^{1}_{A,\Lambda}}$.
This will be a further extention of the works \cite{dz:1,r:1,Sh-T:0}.
For this application, we obtain the following useful proposition which might enable us to derive the contradiction if the coefficient in the condition {\rm (vi)}
were zero:
\begin{prop*}[Proposition~\ref{prop:coeff 0}]
Assume that a Frobenius manifold $M$ of rank $\mu_A$ and dimension one with flat coordinates 
$(t_1,t_{1,1},\dots ,t_{i,j},\dots ,t_{r,a_r-1},t_{\mu_A})$ satisfies the following conditions$:$
\begin{enumerate}
\item 
The unit vector field $e$ and the Euler vector field $E$ are given by
\[
e=\frac{\p}{\p t_1},\ E=t_1\frac{\p}{\p t_1}+\sum_{i=1}^{r}\sum_{j=1}^{a_i-1}\frac{a_i-j}{a_i}t_{i,j}\frac{\p}{\p t_{i,j}}
+\chi_A\frac{\p}{\p t_{\mu_A}}.
\]
\item 
The non--degenerate symmetric bilinear form $\eta$ on $\T_M$ satisfies
\begin{align*}
&\ \eta\left(\frac{\p}{\p t_1}, \frac{\p}{\p t_{\mu_A}}\right)=
\eta\left(\frac{\p}{\p t_{\mu_A}}, \frac{\p}{\p t_1}\right)=1,\\ 
&\ \eta\left(\frac{\p}{\p t_{i_1,j_1}}, \frac{\p}{\p t_{i_2,j_2}}\right)=
\begin{cases}
\frac{1}{a_{i_1}}\quad i_1=i_2\text{ and }j_2=a_{i_1}-j_1,\\
0 \quad \text{otherwise}.
\end{cases}
\end{align*}
\item 
The Frobenius potential $\F$ satisfies $E\F|_{t_{1}=0}=2\F|_{t_{1}=0}$,
\[
\left.\F\right|_{t_1=0}\in\CC\left[[t_{1,1}, \dots, t_{1,a_1-1}, 
\dots, t_{i,j},\dots, t_{r,1}, \dots, t_{r,a_r-1},e^{t_{\mu_A}}]\right].
\]
\item Assume the condition {\rm (iii)}. we have
\begin{equation*}
\F|_{t_1=e^{t_{\mu_A}}=0}=\sum_{i=1}^{r}\G^{(i)}, \quad \G^{(i)}\in \CC[[t_{i,1},\dots, t_{i,a_i-1}]],\ i=1,\dots,r.
\end{equation*}
\item 
Assume the condition {\rm (iii)}. In the frame $\frac{\p}{\p t_1}, \frac{\p}{\p t_{1,1}},\dots, 
\frac{\p}{\p t_{r,a_r-1}},\frac{\p}{\p t_{\mu_A}}$ of $\T_M$,
the product $\circ$ can be extended to the limit $t_1=t_{1,1}=\dots=t_{r,a_r-1}=e^{t_{\mu_A}}=0$.
The $\CC$-algebra obtained in this limit is isomorphic to
\[
\CC[x_1,x_2,\dots, x_r]\left/\left(x_ix_j, \ a_ix_i^{a_i}-a_jx_j^{a_j}
\right)_{1\le i\ne j\le r}\right.,
\]
where $\p/\p t_{i,j}$ are mapped to
$x^{j}_i$ for $i=1,\dots,r, j=1,\dots, a_{i}-1$ and $\p/\p t_{\mu_A}$ are mapped to $a_{1}x_{1}^{a_1}$.
\item The term 
\[
\displaystyle\left(\prod_{i=1}^{r}t_{i,1}\right)e^{t_{\mu_A}}
\]
occurs with the coefficient $0$ in $\F$. 
\item 
The term $t_{i,1}^{2}t_{i,a_{i}-1}^{2}$ in $\F$ occurs with the coefficient 
\[
\begin{cases}
\displaystyle
-1/96 \quad \text{if} \quad a_{i}=2,\\
\displaystyle
-1/4a_{i}^{2} \quad \text{if} \quad a_{i}\ge3.
\end{cases}
\]
\end{enumerate}
Then any term $t^{\alpha}e^{mt_{\mu_A}}$ for $m\ge 1$ occurs with the coefficient $0$ in $\F$.
\end{prop*}
\bigskip
\noindent
{\it Acknowledgement}\\
\indent
The author would express his deep gratitude to Professor Atsushi Takahashi 
for his valuable discussion and encouragement.   


\section{Preliminary}\label{section2}

In this section, we recall the definition and three elementary properties of the Frobenius manifold \cite{d:1}. 
The definition below is taken from Saito--Takahashi \cite{st:1}.
\begin{defn}
Let $M=(M,\O_{M})$ be a connected complex manifold or a formal manifold over $\CC$ of dimension $\mu$
whose holomorphic tangent sheaf and cotangent sheaf are denoted by $\T_{M}$ and $\Omega_M^1$ respectively.
Set a complex number $d$.
A {\it Frobenius structure of
rank $\mu$ and dimension $d$ on M} is a tuple $(\eta, \circ , e,E)$, where we denote by $\eta$ a non--degenerate $\O_{M}$--symmetric bilinear
form on $\T_{M}$, by $\circ $ an $\O_{M}$-bilinear product on $\T_{M}$ of an associative and commutative
$\O_{M}$--algebra structure with the unit $e$ and by $E$ a holomorphic vector field on $M$ called
the Euler vector field, satisfying the following axioms:
\begin{enumerate}
\item The product $\circ$ is self--ajoint with respect to $\eta$: that is,
\begin{equation*}
\eta(\delta\circ\delta',\delta'')=\eta(\delta,\delta'\circ\delta''),\quad
\delta,\delta',\delta''\in\T_M. 
\end{equation*} 
\item The {\rm Levi}--{\rm Civita} connection $\ns:\T_M\otimes_{\O_M}\T_M\to\T_M$ with respect to $\eta$ is
flat: that is, 
\begin{equation*}
[\ns_\delta,\ns_{\delta'}]=\ns_{[\delta,\delta']},\quad \delta,\delta'\in\T_M.
\end{equation*}
\item The tensor $C:\T_M\otimes_{\O_M}\T_M\to \T_M$  defined by 
$C_\delta\delta':=\delta\circ\delta'$, $(\delta,\delta'\in\T_M)$ is flat: that is,
 \begin{equation*}
\ns C=0.
\end{equation*} 
\item The unit $e$ for the product $\circ$ is a $\ns$--flat holomorphic vector field: that is,
\begin{equation*}
\ns e=0.
\end{equation*} 
\item The non--degenerate bilinear form $\eta$ and the product $\circ$ are homogeneous of degree $2-d$ and $1$ respectively 
with respect to the Lie derivative ${\rm Lie}_{E}$ of the {\rm Euler} vector field $E$: that is,
\begin{equation*}
{\rm Lie}_E(\eta)=(2-d)\eta,\quad {\rm Lie}_E(\circ)=\circ.
\end{equation*}
\end{enumerate}
\end{defn}
We shall expose, without their proof, three basic properties of the Frobenius manifold which are necessary to state Theorem \ref{first}.
Let us consider the space of horizontal sections of the connection $\ns$:
\[
\T_M^f:=\{\delta\in\T_M~|~\ns_{\delta'}\delta=0\text{ for all }\delta'\in\T_M\}.
\]
Then the axiom {\rm (ii)} implies that $\T_M^f$ is a local system of rank $\mu $ on $M$: 
\begin{prop}\label{prop:flat coordinates}
At each point of the Frobenius manifold $M$, there exists a local coordinate $(t_1,\dots,t_{\mu})$, called flat coordinates, such that
$e=\p_1$, $\T_M^f$ is spanned by $\p_1,\dots, \p_{\mu}$ and $\eta(\p_i,\p_j)\in\CC$ for all $i,j=1,\dots, \mu$
where we denote $\p/\p t_i$ by $\p_i$. 
\end{prop} 
The axiom {\rm (iii)} implies the existence of the Frobenius potential:
\begin{prop}\label{prop:potential}
At each point of the Frobenius manifold $M$, there exists the local holomorphic function $\F$, called Frobenius potential, satisfying
\begin{equation*}
\eta(\p_i\circ\p_j,\p_k)=\eta(\p_i,\p_j\circ\p_k)=\p_i\p_j\p_k \F,
\quad i,j,k=1,\dots,\mu,
\end{equation*}
for any system of flat coordinates. In particular, we have
\begin{equation*}
\eta_{ij}:=\eta(\p_i,\p_j)=\p_1\p_i\p_j \F. 
\end{equation*}
\end{prop}
Furthermore, the associativity of the product $\circ$ implies that the Frobenius potential satisfies the WDVV equations:
\begin{prop}\label{prop:wdvv-equation}
The Frobenius potential $\F$ satisfies the following equations:
\[
\displaystyle \sum_{\sigma ,\tau=1}^{\mu}\p_{a}\p_{b}\p_{\sigma}\F \cdot \eta^{\sigma \tau}\cdot \p_{\tau}\p_{c}\p_{d}\F
-\sum_{\sigma ,\tau=1}^{\mu}\p_{a}\p_{c}\p_{\sigma}\F \cdot \eta^{\sigma \tau}\cdot \p_{\tau}\p_{b}\p_{d}\F=0,
\]
where $a, b, c, d\in \{1,\dots, \mu\}$.
\end{prop}
\section{A Uniqueness Theorem}
Let $r\ge 3$ be a positive integer.
Let $A$ be a multiplet $(a_1,a_2,\dots, a_r)$ of positive integers such that $2\le a_1\le a_2\le \dots \le a_r$ 
and $\Lambda=(\lambda_1,\ldots,\lambda_r)$ a multiplet of pairwise distinct elements of $\PP^{1}(\CC)$ normalized such that 
$\lambda_1=\infty$, $\lambda_2=0$ and $\lambda_3=1$.
Set
\begin{equation}
\mu_A:=2+\sum_{k=1}^{r} \left(a_{k}-1\right),
\end{equation}
\begin{equation}
\chi_A:=2+\sum_{k=1}^{r}\left(\frac{1}{a_k}-1\right).
\end{equation}
We have the following uniqueness theorem for Frobenius manifolds of rank $\mu_A$ and dimension one. 
The proof of this uniqueness theorem, especially Proposition~\ref{mreconst}, is inspired by Kontsevich--Manin
\cite{km:1} and E. Mann \cite{mann:1}:
\begin{thm}\label{first}
There exists a unique Frobenius manifold $M$ of rank $\mu_A$ and dimension one with flat coordinates 
$(t_1,t_{1,1},\dots ,t_{i,j},\dots ,t_{r,a_r-1},t_{\mu_A})$ satisfying the following conditions$:$
\begin{enumerate}
\item 
The unit vector field $e$ and the Euler vector field $E$ are given by
\[
e=\frac{\p}{\p t_1},\ E=t_1\frac{\p}{\p t_1}+\sum_{i=1}^{r}\sum_{j=1}^{a_i-1}\frac{a_i-j}{a_i}t_{i,j}\frac{\p}{\p t_{i,j}}
+\chi_A\frac{\p}{\p t_{\mu_A}}.
\]
\item 
The non--degenerate symmetric bilinear form $\eta$ on $\T_M$ satisfies
\begin{align*}
&\ \eta\left(\frac{\p}{\p t_1}, \frac{\p}{\p t_{\mu_A}}\right)=
\eta\left(\frac{\p}{\p t_{\mu_A}}, \frac{\p}{\p t_1}\right)=1,\\ 
&\ \eta\left(\frac{\p}{\p t_{i_1,j_1}}, \frac{\p}{\p t_{i_2,j_2}}\right)=
\begin{cases}
\frac{1}{a_{i_1}}\quad i_1=i_2\text{ and }j_2=a_{i_1}-j_1,\\
0 \quad \text{otherwise}.
\end{cases}
\end{align*}
\item 
The Frobenius potential $\F$ satisfies $E\F|_{t_{1}=0}=2\F|_{t_{1}=0}$,
\[
\left.\F_{A}\right|_{t_1=0}\in\CC\left[[t_{1,1}, \dots, t_{1,a_1-1}, 
\dots, t_{i,j},\dots, t_{r,1}, \dots, t_{r,a_r-1},e^{t_{\mu_A}}]\right].
\]
\item Assume the condition {\rm (iii)}. we have
\begin{equation*}
\F_{A}|_{t_1=e^{t_{\mu_A}}=0}=\sum_{i=1}^{r}\G_{A}^{(i)}, \quad \G_{A}^{(i)}\in \CC[t_{i,1},\dots, t_{i,a_i-1}],\ i=1,\dots,r.
\end{equation*}
\item 
Assume the condition {\rm (iii)}. In the frame $\frac{\p}{\p t_1}, \frac{\p}{\p t_{1,1}},\dots, 
\frac{\p}{\p t_{r,a_r-1}},\frac{\p}{\p t_{\mu_A}}$ of $\T_M$,
the product $\circ$ can be extended to the limit $t_1=t_{1,1}=\dots=t_{r,a_r-1}=e^{t_{\mu_A}}=0$.
The $\CC$-algebra obtained in this limit is isomorphic to
\[
\CC[x_1,x_2,\dots, x_r]\left/\left(x_ix_j, \ a_ix_i^{a_i}-a_jx_j^{a_j}
\right)_{1\le i\ne j\le r}\right.,
\]
where $\p/\p t_{i,j}$ are mapped to
$x^{j}_i$ for $i=1,\dots,r, j=1,\dots, a_{i}-1$ and $\p/\p t_{\mu_A}$ are mapped to $a_{1}x_{1}^{a_1}$.
\item The term 
\[
\displaystyle\left(\prod_{i=1}^{r}t_{i,1}\right)e^{t_{\mu_A}}
\]
occurs with the coefficient $1$ in $\F_{A}$. 
\end{enumerate}
\end{thm}
\begin{rem}
The conditions in Theorem~\ref{first} are satisfied by natural ones for
the orbifold Gromow--Witten theory of $\PP^{1}_{A,\Lambda}$. The condition {\rm (i)}, {\rm (ii)} and {\rm (v)} come 
from the conditions for a homogeneous basis of the orbifold cohomology group, the orbifold Poincar\'e pairing
and the large radius limit for the orbifold Gromov--Witten theory respectively. The condition 
{\rm (ii)} and {\rm (v)} are essential to obtain coefficients corresponding to
genus zero three points degree zero correlators. The condition {\rm (iii)} comes from
the divisor axiom. The condition {\rm (iv)} and {\rm (vi)} come from some geometrical meanings
of the orbifold Gromov--Witten invariants.
Namely, the coefficient of the term in the condition {\rm (vi)} corresponds to a certain 
genus zero $r$--points degree one correlator.   
\end{rem}
We shall notice different and common points between the present proof of Theorem \ref{first} and the one for Theorem 3.1 in \cite{ist:1}.
Surprisingly, Theorem \ref{first} can be proven by the parallel way to the one in our previous paper \cite{ist:1}.
However, for general cases $r\ge 4$, we have to modify the arguments in \cite{ist:1} related to the reconstruction of the coefficients
corresponding to the genus zero degree one correlators, e.g., the terms in the WDVV equations whose coefficients give the recursion relations.
In particular, the arguments in Proposition~\ref{0reconst} except for Lemma~\ref{lem4-0reconst}
are very natural generalizations of the one for Proposition 3.36 in \cite{ist:1}. 
In contrast to this, some arguments in \cite{ist:1} can be also applied without major modifications, e.g., the arguments in \cite{ist:1} 
related to the reconstruction of the coefficients corresponding to the genus zero higher degree correlators.
From now on, we will mark with asterisks ($*$) on propositions, lemmas and sublemmas whose proofs are (almost) same with the ones in \cite{ist:1}
and mark with daggers ($\dagger$) on them whose proofs need some modifications.
In order to make the proof self--contained, we shall include all details of arguments even if the arguments are common to the ones in 
our previous paper \cite{ist:1}. 

We shall use the same notations with the ones in our previous paper \cite{ist:1}.
By the condition {\rm (iii)} in Theorem~\ref{first},
we can expand the non--trivial part of the Frobenius potential $\F_{A}|_{t_1=0}$ as
\[
\F_{A}|_{t_{1}=0}=\sum_{\alpha =(\alpha_{1,1},\dots ,\alpha_{r,a_{r}-1})} c(\alpha ,m) t^{\alpha}e^{mt_{\mu_A}}, \ \
t^{\alpha}=\prod_{i=1}^{r}\prod_{j=1}^{a_i-1}t_{i,j}^{\alpha_{i,j}}.
\]
Here we note that, by Proposition~\ref{prop:potential}, 
the terms in $\F_A$ including $t_{1}$ are only cubic terms $t_{1}t_{i,j}t_{i,a_{i}-j}$ and their coefficients can be determined by the condition {\rm (ii)}.

Consider a free abelian group $\ZZ^{\mu_A-2}$ and denote  its 
standard basis by $e_{i,j}$, $i=1,\dots, r$, $j=1,\dots, a_i-1$.
The element $\alpha=\sum_{i,j}\alpha_{i,j}e_{i,j}$, $\alpha_{i,j}\in\ZZ$ of  $\ZZ^{\mu_A-2}$ is called {\it non--negative} 
and is denoted by  $\alpha\ge 0$ if all $\alpha_{i,j}$ are non--negative integers.
We also denote by $c(e_{1}+e_{i,j}+e_{i,a_{i}-j},0)$ the coefficient of $t_{1}t_{i,j}t_{i,a_{i}-j}$ in the trivial part of the Frobenius potential $\F$.
For a non--negative $\alpha\in\ZZ^{\mu_A-2}$, we set 
\[
|\alpha|:=\sum_{i=1}^{r}\sum_{j=1}^{a_i-1}\alpha_{i,j},
\]
and call it the {\it length} of $\alpha$. 
Define the number $s_{a,b,c}$ for $a,b,c\in \ZZ$ as follows:
\[
s_{a,b,c}=
\begin{cases}
1 \ \ \text{if} \ \ a,b,c \text{ are pairwise distinct},\\
6 \ \ \text{if} \ \ a=b=c,\\
2 \ \ \text{otherwise}.\\
\end{cases}
\]
For $a, b, c, d\in \{1,\dots, \mu_A\}$, denote by $WDVV(a,b,c,d)$   
the following equation:
\[
\displaystyle \sum_{\sigma ,\tau=1}^{\mu_A}\p_{a}\p_{b}\p_{\sigma}\F \cdot \eta^{\sigma \tau}\cdot \p_{\tau}\p_{c}\p_{d}\F
-\sum_{\sigma ,\tau=1}^{\mu_A}\p_{a}\p_{c}\p_{\sigma}\F \cdot \eta^{\sigma \tau}\cdot \p_{\tau}\p_{b}\p_{d}\F=0,
\]
where $(\eta^{\sigma \tau}):=(\eta_{\sigma \tau})^{-1}$.

\subsection{Coefficients $c(\alpha,0)$ and $c(\alpha,1)$ can be reconstructed}\label{reconst-s1}
\begin{propa}\label{lem3}
Coefficients $c(\alpha,0)$ with $|\alpha|=3$ are determined by the condition {\rm (v)} of Theorem~\ref{first}.
\end{propa}
\begin{pf}
Note that $C_{ijk}=\eta (\partial_{i}\circ \partial_{j}, \partial_{k})$
and the non--degenerate bilinear form $\eta$ can be extended to the limit $\underline{t},e^{t} \rightarrow 0$.
We denote by $\eta'$ this extended bilinear form. 
By the condition {\rm (v)}, the relation $x_{i}x_{j}=0$ if $i\ne j$ holds in the $\CC$--algebra obtained in this limit.
Therefore, we have $c(\sum_{k=1}^3e_{i_k,j_{k}},0)\ne 0$ only if $i_1=i_2=i_3$.
In particular, we have 
\begin{multline*}
s_{j_{1},j_{2},j_{3}}\cdot c\left(\sum_{k=1}^3e_{i,j_{k}},0\right)=\lim_{\underline{t},e^{t} \rightarrow 0} \p_{i,j_{1}}\p_{i,j_{2}}\p_{i,j_{3}}\F_{A}
=\eta' (x_{i}^{j_{1}}\cdot x_{i}^{j_{2}},x_{i}^{j_{3}})\\
=\eta' (1\cdot x_{i}^{j_1+j_2}, x_{i}^{j_{3}})
=\lim_{\underline{t},e^{t} \rightarrow 0} \p_{1}\p_{i, j_{1}+j_{2}}\p_{i,j_{3}}\F_{A}
\end{multline*}
by Proposition~\ref{prop:potential} and 
\[
\lim_{\underline{t},e^{t} \rightarrow 0} \p_{1}\p_{i, j_{1}+j_{2}}\p_{i, j_{3}}\F=
\begin{cases}
\frac{1}{a_i}\quad \text{if } \sum_{k=1}^3j_{k}=a_i,\\
0\quad \text{otherwise}.
\end{cases}
\]
\qed
\end{pf}

\begin{propa}\label{lem3.1-ver1}
A coefficient $c(\alpha,1)$ with $|\alpha|\le r$ 
is none--zero if and only if $\alpha=\sum_{k=1}^{r}e_{k,1}$.
In particular, we have $c(\sum_{k=1}^{r}e_{k,1},1)=1$  
by the condition {\rm (vi)} of Theorem~\ref{first}.
\end{propa}
\begin{pf}

We shall split the proof into following two cases.

\begin{lema}[Case 1]\label{lem1-lem3.1}
Let $\gamma \in \ZZ^{\mu_{A}-2}$ be a non--negative element satisfying that $|\gamma|=r$ and $\gamma-e_{i,j}\geq 0$ for some $i, j$.
If $a_{i}\ge 3$ and  $j\ge 2$, then we have $c(\gamma,1)=0$.
\end{lema}
\begin{pf}
Since $\deg(t^{\alpha}e^{t_{\mu_{A}}})<2$, we have $c(\alpha, 1)=0$ if $|\alpha|\le r-1$.
We shall calculate the coefficient of the term  
$t^{\gamma-e_{i,j}}e^{t_{\mu_A}}$
in $WDVV((i,1),(i,j-1),\mu_A,\mu_A)$. Then we have
\[
s_{1,j-1,a_{i}-j}\cdot c(e_{i,1}+e_{i,j-1}+e_{i,a_{i}-j},0)\cdot a_{i} \cdot \gamma_{i,j} \cdot c(\gamma ,1)=0.
\]
Hence we have $c(\gamma ,1)=0$.
\qed
\end{pf}

\begin{lema}[Case 2]\label{lem2-lem3.1}
If a non--negative element $\gamma \in \ZZ^{\mu_{A}-2}$ satisfies that
$|\gamma|=r$ and $\gamma =\sum_{k=1}^{r}\gamma_{k,1}e_{k,1}$ 
for some $\gamma_{1,1}, \dots, \gamma_{r,1}$ such that $\prod_{k=1}^{r}\gamma_{k,1}=0$,
then we have $c(\gamma ,1)=0$.
\end{lema}
\begin{pf}
Note that $c(\alpha,0)=0$ if $|\alpha|=4$ and $\alpha-e_{i_{1},j_{1}}-e_{i_{2},j_{2}}\ge 0$ for $i_{1}\neq i_{2}$
by the condition {\rm (iv)} and that $c(\alpha, 1)=0$ if $|\alpha|\le r-1$
since $\deg(t^{\alpha}e^{t_{\mu_{A}}})<2$. Assume that $\gamma_{i,1}=0$.
We shall calculate the coefficient of the term $(\prod_{k\ne i}^{r}t_{k,1}^{\gamma_{k,1}})e^{t_{\mu_A}}$
in the WDVV equation $WDVV((i,1),(i,a_{i}-1),\mu_A,\mu_A)$. Then we have
\[
c(e_{1}+e_{i,a_{i}-1}+e_{i,1},0)\cdot c(\gamma ,1)=0.
\]
Hence we have $c(\gamma ,1)=0$ and hence Lemma~\ref{lem2-lem3.1}.
\qed
\end{pf}
Therefore we have Proposition~\ref{lem3.1}.
\qed
\end{pf}

\begin{cora}\label{lem4-1}
If $a_{i}\ge 3$, then we have 
\[
c(2e_{i,1}+2e_{i,a_{i}-1},0)=-\frac{1}{4a_{i}^{2}}. 
\]
\end{cora}
\begin{pf}
By the condition {\rm (iv)},
we have $c(\gamma,0)=0$ if $\gamma-e_{i_1,j_1}-e_{i_2,j_2}\ge 0$ for $i_1\ne i_2$. 
We shall calculate the coefficient of the term $(\prod^{r}_{k=1} t_{k,1})e^{t_{\mu_A}}$ in $WDVV((i,1),(i,a_{i}-1),\mu_A ,\mu_A)$. 
Then we have
\begin{multline*}
c(e_{1}+e_{i,1}+e_{i,a_{i}-1},0)\cdot 1\cdot c(\sum_{k=1}^{r} e_{k,1},1)+\\
4\cdot c(2e_{i,1}+2e_{i,a_{i}-1},0)\cdot a_{i}\cdot c(\sum_{k=1}^{r} e_{k,1},1)=0.
\end{multline*}
We have $c(e_{1}+e_{i,1}+e_{i,a_{i}-1},0)=1/a_{i}$ and
$c(\sum_{k=1}^{r} e_{k,1},1)=1$ by the conditions {\rm (ii)} and {\rm (vi)} in Theorem~\ref{first}.
Hence we have $c(2e_{i,1}+2e_{i,a_{i}-1},0)=-1/4a_{i}^{2}$.
\qed
\end{pf}

\begin{cora}\label{lem4-2}
If $a_{i}=2$,
then we have  
\[
c(4e_{i,1},0)=-\frac{1}{96}. 
\]
\end{cora}
\begin{pf}
By the condition {\rm (iv)}, we have $c(\gamma,0)=0$ if $\gamma-e_{i_1,j_1}-e_{i_2,j_2}\ge 0$ for $i_1\ne i_2$.
We shall calculate the coefficient of the term $(\prod^{r}_{k=1} t_{k,1})e^{t_{\mu_A}}$ in $WDVV((i,1),(i,1),\mu_A ,\mu_A)$.
Then we have
\[
2c(e_{1}+2e_{i,1},0)\cdot c(\sum_{k=1}^{r} e_{k,1},1)+24c(2e_{i,1}+2e_{i,a_{i}-1},0)\cdot 2\cdot c(\sum_{k=1}^{r} e_{k,1},1)=0.
\]
We have $c(e_{1}+2e_{i,1},0)=1/4$ and $c(\sum_{k=1}^{r} e_{i,1},1)=1$ by the conditions {\rm (ii)} and {\rm (vi)} in Theorem~\ref{first}.
Hence we have $c(4e_{i,1},0)=-1/96$.
\qed
\end{pf}

\begin{propd}\label{0reconst}
Assume that $c(\alpha,0)$ and $c(\alpha',1)$ are reconstructed if $|\alpha|\le k+3$ and $|\alpha'|\le k+r$ for some $k\in\ZZ_{\ge 0}$.
Then coefficients $c(\gamma,0)$ with $|\gamma |\le k+4$ and $c(\gamma',1)$ with $|\gamma' |\le k+r+1$ are reconstructed
from coefficients $c(\alpha,0)$ with $|\alpha|\le k+3$ and $c(\alpha',1)$ with $|\alpha'|\le k+r$.
\end{propd}
\begin{pf}
We shall split the proof of Proposition \ref{0reconst} into following four steps. 

\begin{lemd}[Step 1]\label{lem1-0reconst}
If a non--negative element $\beta \in \ZZ^{\mu_{A}-2}$ satisfies that $|\beta|=k+1$,
then the coefficient $c(\beta +e_{i,j}+e_{i,j'}+e_{i,a_{i}-1},0)$ for some $i,j,j'$ can be reconstructed from coefficients
$c(\alpha,0)$ with $|\alpha|\le k+3$ and $c(\alpha',1)$ with $|\alpha'|\le k+r$.
\end{lemd}

\begin{pf}
Without loss of generality, we can assume $i=1$.
First we shall show that the coefficient $c(\beta +e_{1,1}+e_{1,j+j'-1}+e_{1,a_{1}-1},0)$ can be reconstructed 
from coefficients $c(\alpha,0)$ with $|\alpha|\le k+3$ and $c(\alpha',1)$ with $|\alpha'|\le k+r$. 
We have $\deg (t^{\beta}t_{1,j+j'-1})=1$. By the condition {\rm (iv)}, there exist $e_{1,l},e_{1,l'}$ such that
\begin{itemize}
\item $\beta +e_{1,j+j'-1}-e_{1,l}-e_{1,l'}\geq 0$, 
\item $\deg(t_{1,l})+\deg(t_{1,l'})\leq 1$.
\end{itemize}
We put $\beta' :=\beta +e_{1,1}+e_{1,j+j'-1}-e_{1,l}-e_{1,l'}$.
We shall calculate the coefficient of the term $t^{\beta'}(\prod_{k=2}^{r}t_{k,1})e^{t_{\mu_A}}$ in the WDVV equation $WDVV((1,l),(1,l'),\mu_A,\mu_A)$.
Then we have
\begin{eqnarray*}
(\beta'_{1,l}+1)(\beta'_{1,l'}+1)(\beta'_{1,a_{1}-1}+1)\cdot c(\beta +e_{1,1}+e_{1,j+j'-1}+e_{1,a_{1}-1},0)
\cdot a_{1}\cdot c(\sum_{k=1}^{r}e_{k,1},1)\\
+(known \ \ terms)=0.
\end{eqnarray*}
By the condition {\rm (vi)} in Theorem~\ref{first}, the coefficient
$c(\beta +e_{1,1}+e_{1,j+j'-1}+e_{1,a_{1}-1},0)$ can be reconstructed from coefficients 
$c(\alpha,0)$ with $|\alpha|\le k+3$ and $c(\alpha',1)$ with $|\alpha'|\le k+r$.

Next we shall show that the coefficient $c(\beta +(\sum_{k=2}^{r}e_{k,1})+e_{1,j+j'},1)$ can be reconstructed  
from coefficients $c(\alpha,0)$ with $|\alpha|\le k+3$ and $c(\alpha',1)$ with $|\alpha'|\le k+r$.
We shall calculate the coefficient of the term  $t^{\beta }(\prod_{k=4}^{r}t_{k,1})e^{t_{\mu_A }}$ in 
the WDVV equation $WDVV((1,1),(1,j+j'-1),(2,1),(3,1))$. 
Then we have 
\begin{eqnarray*}
\lefteqn{s_{1,j+j'-1,a_{1}-j-j'}\cdot c(e_{1,1}+e_{1,j+j'-1}+e_{1,a_{1}-j-j'},0)\cdot a_{1}\cdot }\\
&&(\beta_{1,j+j'}+1)(\beta_{2,1}+1)(\beta_{3,1}+1)\cdot 
c(\beta +(\sum_{k=2}^{r}e_{k,1})+e_{1,j+j'},1)\\
&&+ (known \ \ terms)\\
&&+(\beta_{1,1}+1)(\beta_{1,j+j'-1}+1)(\beta_{1,a_{1}-1}+1)\cdot c(\beta +e_{1,1}+e_{1,j+j'-1}+e_{1,a_{1}-1},0)
\cdot a_{1}\cdot \\
&&c(\sum_{k=1}^{r}e_{k,1},1)=0.
\end{eqnarray*}
By the previous argument and Proposition \ref{lem3},  the coefficient $c(\beta +(\sum_{k=2}^{r}e_{k,1})+e_{1,j+j'},1)$ can be reconstructed 
from $c(\alpha,0)$ with $|\alpha|\le k+3$ and $c(\alpha',1)$ with $|\alpha'|\le k+r$.

Finally we shall show that the coefficient $c(\beta +e_{1,j}+e_{1,j'}+e_{1,a_{1}-1},0)$ can be reconstructed  
from coefficients $c(\alpha,0)$ with $|\alpha|\le k+3$ and $c(\alpha',1)$ with $|\alpha'|\le k+r$.
We shall calculate the coefficient of the term  $t^{\beta }(\prod_{k=2}^{r}t_{k,1})e^{t_{\mu_A }}$ in $WDVV((1,j),(1,j'),\mu_A,\mu_A)$. 
Then we have
\begin{eqnarray*}
\lefteqn{\hspace{-33mm}{\rm (i)} \ (\beta_{1,j}+1)(\beta_{1,j'}+1)(\beta_{1,a_{1}-1}+1)\cdot
c(\beta +e_{1,j}+e_{1,j'}+e_{1,a_{1}-1},0)\cdot a_{1}\cdot c(\sum_{k=1}^{r}e_{k,1},1)}\\ 
&&+ (known \ \ terms)\\ 
&&+s_{j,j'',a_{1}-j-j'}\cdot c(e_{1,j}+e_{1,j'}+e_{1,a_{1}-j-j'},0)\cdot \\
&&a_{1}\cdot (\beta_{1,j+j'}+1)\cdot c(\beta +(\sum_{k=2}^{r}e_{k,1})+e_{1,j+j'},1)=0\\ 
&&\text{if} \ \ a_{1}-j+a_{1}-j'\geq a_{1}+1, 
\end{eqnarray*}
\begin{eqnarray*}
\lefteqn{\hspace{-33mm}{\rm (ii)} \ (\beta_{1,j}+1)(\beta_{1,j'}+1)(\beta_{1,a_{1}-1}+1)\cdot
c(\beta +e_{1,j}+e_{1,j'}+e_{1,a_{1}-1},0)\cdot a_{1}\cdot c(\sum_{k=1}^{r}e_{k,1},1)}\\ 
&&+ (known \ \ terms)\\ 
&&+c(e_{1,j}+e_{1,j'}+e_{1},0)\cdot 1\cdot c(\beta +(\sum_{k=2}^{r}e_{k,1}),1)=0\\
&&\text{if} \ \ a_{1}-j+a_{1}-j'=a_{1}, \\
\lefteqn{\hspace{-33mm}{\rm (iii)} \ (\beta_{1,j}+1)(\beta_{1,j'}+1)(\beta_{1,a_{1}-1}+1)\cdot
c(\beta +e_{1,j}+e_{1,j'}+e_{1,a_{1}-1},0)\cdot a_{1}\cdot c(\sum_{k=1}^{r}e_{k,1},1)}\\ 
&&+ (known \ \ terms)=0\\ 
&&\text{if} \ \ a_{1}-j+a_{1}-j'<a_{1}.
\end{eqnarray*}
By the second argument and Proposition \ref{lem3}, the coefficient
$c(\beta +e_{1,j}+e_{1,j'}+e_{1,a_{1}-1},0)$ can be reconstructed from coefficients $c(\alpha,0)$ with $|\alpha|\le k+3$ and $c(\alpha',1)$ with $|\alpha'|\le k+r$.
\qed
\end{pf}

\begin{lemd}[Step 2]\label{lem2-0reconst}
For a non--negative $\gamma \in \ZZ^{\mu_{A}-2}$ with $|\gamma|=k+r+1$,
a coefficient $c(\gamma,1)$ can be reconstructed 
from coefficients $c(\alpha,0)$ with $|\alpha|\le k+3$ and $c(\alpha',1)$ with $|\alpha'|\le k+r$.
\end{lemd}
\begin{pf}
We shall split the proof of Lemma \ref{lem2-0reconst} into following three cases.  
\begin{sublemd}[Step 2--Case 1]\label{sl1-l2-0reconst}
If a non--negative element $\gamma \in \ZZ^{\mu_{A}-2}$ satisfies that $|\gamma|=k+r+1$ and
$\gamma-e_{i,j}\geq 0$ for some $i,j$ such that $j\ge 2$,
then the coefficient $c(\gamma,1)$ can be reconstructed from coefficients $c(\alpha,0)$ with $|\alpha|\le k+3$ and $c(\alpha',1)$ with $|\alpha'|\le k+r$.
\end{sublemd}
\begin{pf}
Put $\gamma':=\gamma-e_{i.j}-(\sum_{k\ne i}^{r}e_{k,1})+e_{i,1}+e_{i,j-1}+e_{i,a_{i}-1}$.
We shall calculate the coefficient of the term  
$t^{\gamma-e_{i,j}}e^{t_{\mu_A}}$
in $WDVV((i,1),(i,j-1),\mu_A,\mu_A)$. Then we have
\begin{eqnarray*}
\lefteqn{s_{1,j-1,a_{i}-j}\cdot c(e_{i,1}+e_{i,j-1}+e_{i,a_{i}-j},0)\cdot a_{i} \cdot 
\gamma_{i_{1},j_{1}} \cdot c(\gamma ,1)}\\ 
&&+ (known\ \ terms)\\
&&+ \gamma'_{i,1}\gamma'_{i,j-1}\gamma'_{i,a_{i}-1}\cdot c(\gamma',0)\cdot a_{i}\cdot c(\sum_{k=1}^{r}e_{k,1},1)=0
\end{eqnarray*}
By Lemma \ref{lem1-0reconst}, the coefficient
$c(\gamma',0)$ can be reconstructed from $c(\alpha,0)$ and $c(\alpha',1)$ with $|\alpha|\le k+3$ and $|\alpha'|\le k+r$.
Hence the coefficient $c(\gamma,1)$ can be reconstructed from
coefficients $c(\alpha,0)$ with $|\alpha|\le k+3$ and $c(\alpha',1)$ with $|\alpha'|\le k+r$.
\qed
\end{pf}
\begin{sublema}[Step 2--Case 2]\label{sl2-l2-0reconst}
If a non--negative element $\gamma \in \ZZ^{\mu_{A}-2}$ satisfies that $|\gamma|=k+r+1$ and
$\gamma=\sum_{k=1}^{r}\gamma_{k,1}e_{k,1}$ for  some  $\gamma_{1,1},\dots, \gamma_{r,1}$ such that $\prod_{k=1}^{r}\gamma_{k,1}\ne 0$,
then we have $c(\gamma,1)=0$.
\end{sublema}
\begin{pf}
By counting the degree of the term $t^{\gamma}e^{t_{\mu_{A}}}$, we have
\[
\deg(t^{\gamma}e^{t_{\mu_{A}}})>\deg((\prod_{k=1}^{r}t_{k,1})e^{t_{\mu_{A}}})=2.
\]
Then we have $c(\gamma,1)=0$.
\qed
\end{pf}
\begin{sublemd}[Step2--Case 3]\label{sl3-l2-0reconst}
If a non--negative element $\gamma \in \ZZ^{\mu_{A-2}}$ satisfies that $|\gamma|=k+r+1$ and $\gamma =\sum_{k=1}^{r}\gamma_{k,1}e_{k,1}$ 
for some  $\gamma_{1,1},\dots,\gamma_{r,1}$ such that  $\prod_{k=1}^{r}\gamma_{k,1}=0$, 
then the coefficient $c(\gamma,1)$ can be reconstructed from coefficients $c(\alpha,0)$ with $|\alpha|\le k+3$ and $c(\alpha',1)$ with $|\alpha'|\le k+r$. 
\end{sublemd}
\begin{pf}
Assume that $\gamma_{i,1}=0$ and put $\gamma':=\gamma-(\sum_{k\ne i}^{r}e_{k,1})+e_{i,1}+e_{i,a_{i}-1}+e_{i,a_{i}-1}$.
We shall calculate the coefficient of the term $(\prod_{k\ne i}^{r}t_{k,1}^{\gamma_{k,1}})e^{t_{\mu_A}}$
in the WDVV equation $WDVV((i,1),(i,a_{i}-1),\mu_A,\mu_A)$. Then we have 
\begin{eqnarray*}
\lefteqn{c(e_{i,1}+e_{1,a_{i}-1}+e_{1},0)\cdot c(\gamma ,1)
+(known \ \ terms)}\\
&&+\gamma'_{i,1}\gamma'_{i,a_{i}-1}\gamma'_{i,a_{i}-1}\cdot c(\gamma',0)\cdot a_{i}\cdot c(\sum_{k=1}^{r}e_{k,1},1)=0
\end{eqnarray*}
We have $\gamma'-e_{k,1}\ge 0$ for some $k\ne i$.
Then we have $c(\gamma',0)=0$ by the condition {\rm (iv)}.
Hence the coefficient $c(\gamma ,1)$ can be reconstructed from coefficients 
$c(\alpha,0)$ with $|\alpha|\le k+3$ and $c(\alpha',1)$ with $|\alpha'|\le k+r$.
Therefore a coefficient $c(\gamma,1)$ can be reconstructed from coefficients
$c(\alpha,0)$ with $|\alpha|\le k+3$ and $c(\alpha',1)$ with $|\alpha'|\le k+r$.
\qed
\end{pf}
Then we have Lemma~\ref{lem2-0reconst}.
\qed
\end{pf}

\begin{lemd}\label{lem4.1}
If $a_{i}\geq 3$, 
then we have  
\[
c(e_{i,j+1}+e_{i,a_{i}-j}+\sum_{k\ne i}^{r}e_{k,1},1)=
\begin{cases}
\frac{1}{a_{i}} \ \ \text{if} \ \ a_{i}-j\neq j+1,\\
\frac{1}{2a_{i}} \ \ \text{if} \ \ a_{i}-j=j+1.
\end{cases}
\]
\end{lemd}
\begin{pf}
We shall calculate the coefficient of the term $(\prod_{k\ne i}^{r} t_{k,1})e^{t_{\mu_A}}$ 
in the WDVV equation $WDVV((i,a_{i}-j-1),(i,1),(i,j+1),\mu_A)$. 
Then we have
\begin{eqnarray*}
\lefteqn{\hspace{-30mm}{\rm (i)} \ \frac{1}{a_{i}}\cdot a_{i} \cdot c(e_{i,j+1}+e_{i,a_{i}-j}+\sum_{k\ne i}^{r}e_{k,1},1)}\\
&&-1\cdot 1\cdot c(e_{1}+e_{i,a_{i}-j-1}+e_{i,j+1},0)=0\\
&&\text{if} \ \ a_{i}-j\neq j+1,\ a_{i}-j-1\neq j+1,
\end{eqnarray*}
\begin{eqnarray*}
\lefteqn{\hspace{-30mm}{\rm (ii)} \ \frac{1}{a_{i}}\cdot a_{i} \cdot c(e_{i,j+1}+e_{i,a_{i}-j}+\sum_{k\ne i_1}^{r}e_{k,1},1)}\\
&&-1\cdot 1\cdot 2\cdot c(e_{1}+e_{i,a_{i}-j-1}+e_{i,j+1},0)=0\\
&&\text{if} \ \ a_{i}-j\neq j+1, \ a_{i}-j-1=j+1,\\
\lefteqn{\hspace{-30mm}{\rm (iii)} \ \frac{1}{a_{i}}\cdot a_{i} \cdot 2\cdot c(e_{i,j+1}+e_{i,a_{i}-j}+\sum_{k\ne i_1}^{r}e_{k,1},1)}\\
&&-1\cdot 1\cdot c(e_{1}+e_{i,a_{i}-j-1}+e_{i,j+1},0)=0\\
&&\text{if} \ \ a_{i}-j=j+1.
\end{eqnarray*}
Hence we have Lemma \ref{lem4.1}.
\qed
\end{pf}

\begin{lemd}[Step 3]\label{lem3-0reconst}
If a non--negative element $\gamma \in \ZZ^{\mu_{A}-2}$ satisfies that 
$|\gamma|=k+4$ and $\gamma-e_{i,1}\geq 0$ for some $i$,
then the coefficient $c(\gamma, 0)$ can be reconstructed from coefficients $c(\alpha,0)$ with $|\alpha|\le k+3$ and $c(\alpha',1)$ with $|\alpha'|\le k+r$.
\end{lemd}
\begin{pf}
We will show this claim by the induction on the degree of parameter $t_{i,j}$.
By Lemma \ref{lem1-0reconst}, the coefficient $c(\beta +e_{i,j}+e_{i,j'}+e_{i,a_{i}-1},0)$ 
with $|\beta|=k+1$ can be reconstructed from coefficients $c(\alpha,0)$ with $|\alpha|\le k+3$ and $c(\alpha',1)$ with $|\alpha'|\le k+r$.
Assume that $c(\gamma',0)$ with $|\gamma'|=k+4$ is known if $\gamma'-e_{i,1}-e_{i,n}\geq 0$, $n\ge l$.

We shall show that a coefficient $c(\gamma,0)$ can be reconstructed from coefficients $c(\alpha,0)$ with $|\alpha|\le k+3$ and $c(\alpha',1)$ with $|\alpha'|\le k+r$
if $|\gamma|=k+4$ and $\gamma-e_{i,1}-e_{i,l-1}\geq 0$.
We have $\deg(t^{\gamma-e_{i,1}-e_{i,l-1}})=l/a_{i}$. 
By the condition {\rm (iv)}, there exist $e_{i,m},e_{i,m'}$ such that 
\begin{itemize}
\item $\gamma-e_{i,1}-e_{i,l-1}-e_{i,m}-e_{i,m'}\geq 0$,
\item $\deg(t_{i,m})+\deg(t_{i,m'})\le l/a_{i}$.
\end{itemize}
Note that $\deg(t_{i,m})+\deg(t_{i,l})<1$ and a coefficient $c(\alpha,1)$ with $|\alpha|=k+r+1$
can be reconstructed by Lemma~\ref{lem2-0reconst}.
We put $\beta:=\gamma-e_{i,l-1}-e_{i,m}-e_{i,m'}$.
We shall calculate the coefficient of the term $t^{\beta}(\prod_{k\ne i}^{r} t_{k,1})e^{t_{\mu_A}}$ in $WDVV((i,m),(i,m'),(i_1,l),\mu_A)$. 
Then we have
\[
\gamma_{i,m}\gamma_{i,m'}\gamma_{i,l-1}\cdot c(\gamma ,0)\cdot a_{i}\cdot c(e_{i,a_{i}+1-l}+e_{i,l}+\sum_{k\ne i}^{r}e_{k,1},1)
+(known \ \ terms)=0.
\]
By Lemma \ref{lem4.1}, the coefficient $c(\gamma ,0)$ can be reconstructed from coefficients $c(\alpha,0)$ with $|\alpha|\le k+3$ and $c(\alpha',1)$ with $|\alpha'|\le k+r$
if $\gamma-e_{i,1}\ge 0$. 
\qed
\end{pf}

\begin{lema}[Step 4]\label{lem4-0reconst}
A coefficient $c(\gamma,0)$ with $|\gamma|=k+4$ can be reconstructed from $c(\alpha,0)$ with $|\alpha|\le k+3$ and $c(\alpha',1)$ with $|\alpha'|\le k+r$.
\end{lema}
\begin{pf}
We will show this claim by the induction on the degree of parameter $t_{i,j}$.
By Lemma \ref{lem3-0reconst}, a coefficient $c(\gamma ,0)$ can be reconstructed from 
$c(\alpha,0)$ with $|\alpha|\le k+3$ and $c(\alpha',1)$ with $|\alpha'|\le k+r$ if $\gamma-e_{i,1}\ge 0$. 
Assume that $c(\gamma',0)$ with $|\gamma'|=k+4$ is known if $\gamma'-e_{i,n}\geq 0$ for $n\leq l$.
We shall show a coefficient $c(\gamma,0)$ can be reconstructed from coefficients $c(\alpha,0)$ with $|\alpha|\le k+3$ and $c(\alpha',1)$ with $|\alpha'|\le k+r$
if $\gamma-e_{i,l+1}\ge 0$.
We shall calculate the coefficient of the term $t^{\gamma-e_{i,j}-e_{i,j'}-e_{i,l+1}}$ in $WDVV((i,1),(i,l),(i,j),(i,j'))$. Then we have 
\[
s_{1,l,a_{i}-1-l}\cdot c(e_{i,1}+e_{i,l}+e_{i,a_{i}-1-l},0)\cdot a_{i}\cdot \gamma_{i,j}\gamma_{i,j'}\gamma_{i,l+1}\cdot c(\gamma ,0)
+(known \ \ terms)=0.
\]
Hence a coefficient $c(\gamma ,0)$ can be reconstructed from coefficients $c(\alpha,0)$ with $|\alpha|\le k+3$ and $c(\alpha',1)$ with $|\alpha'|\le k+r$.
\qed
\end{pf}
Therefore we have Proposition~\ref{0reconst}
\qed
\end{pf}
By Proposition \ref{lem3}, Proposition \ref{lem3.1} and Proposition \ref{0reconst}, 
coefficients $c(\gamma,0)$ and $c(\gamma,1)$ can be reconstructed from $c(\beta,0)$ with $|\beta|=3$. 

Let $\F_{A_{i}}(t'_{1}, {\bf t'}_{3}, e^{t_{\mu_{A_{i}}}})$ be the Frobenius potential for the Frobenius manifold $M_{A_{i}}$
in Theorem \ref{first} in \cite{ist:1} where a multiplet of positive integers $A_{i}$ is $(1,1,a_{i})$ such that $a_{i}\ge 2$
and we denote by $(t'_{1}, {\bf t'}_{3}, t_{\mu_{A_{i}}}):=(t'_{1},t'_{3,1},\dots,t'_{3,a_{i}-1}, t_{\mu_{A_{i}}})$ 
the flat coordinate for the Frobenius manifold $M_{A_{i}}$.  
Inspired by Proposition~\ref{lem3}, Corollary~\ref{lem4-1}, Corollary~\ref{lem4-2} and Lemma~\ref{lem4.1}, we have the following Proposition~\ref{atype decomp-prop}.
\begin{prop}\label{atype decomp-prop}
For the polynomial $\G_{A}^{(i)}$ in the condition {\rm (iv)} in Theorem \ref{first}, we have 
\begin{equation*}
\G_{A}^{(i)}=\F_{A_{i}}(0,{\bf t}_{i},0)
\end{equation*}
where ${\bf t}_{i}:=(t_{i,1},\dots,t_{i,a_{i}-1})$ is the $i$--th parts of the flat coordinate in Theorem~\ref{first}.
\end{prop}

\begin{pf}
We can expand the Frobenius potential $\F_{A_{i}}(0, {\bf t'}_{3}, e^{t_{\mu_{A_{i}}}})$ uniquely as follows:
\[
\F_{A_{i}}(0, {\bf t'}_{3}, e^{t_{\mu_{A_{i}}}})=\sum_{\alpha =(\alpha_{i,1},\dots ,\alpha_{i,a_{i}-1})} c'(\alpha ,m) {t'}^{\alpha}e^{mt_{\mu_{A_{i}}}}, \ \
{t'}^{\alpha}=\prod_{j=1}^{a_i-1}{t'}_{3,j}^{\alpha_{i,j}}.
\]
The coefficients $c'(\alpha ,m)$ are uniquely determined by Theorem 3.1 in \cite{ist:1}.
We already proved $c(\alpha, 0)=c'(\alpha,0)$ if $|\alpha|=3$ and $c(\alpha'+\sum^{r}_{k\ne i} e_{k,1},1)=c'(\alpha',1)$ if $|\alpha|=1$
in Proposition~\ref{lem3} and Proposition~\ref{lem3.1-ver1}.
We shall show the proposition by the induction concerning the length and split the proof into the following four steps.
\begin{lem}[Step 1]\label{atype decomp-1}
Assume that $c(\alpha,0)=c'(\alpha,0)$ and $c(\alpha'+\sum^{r}_{k\ne i} e_{k,1},1)=c'(\alpha',1)$
if $|\alpha|\le k+2$ and $|\alpha'|\le k$ for some $k\in \NN$.
Then we have $c(\beta+\sum^{r}_{k\ne i} e_{k,1},1)=c'(\beta,1)$ if $|\beta|=k+1$.
\end{lem}
\begin{pf}
If $\beta-e_{i,1}\ge 0$, then we have $c(\beta+\sum^{r}_{k\ne i} e_{k,1},1)=c'(\beta,1)=0$ since  
both ${\rm deg}(t^{\beta}\prod_{k\ne i}t_{k,1}e^{t_{\mu_A}})$ and ${\rm deg}(t^{\beta}e^{t_{\mu_{A_{i}}}})$ 
are greater than $2$. Hence we have $\beta-e_{i,j}\ge 0$ for some $j\ge 2$.
The coefficient of the term $t^{\beta-e_{i,j}}(\prod_{k\ne i}t_{k,1})e^{t_{\mu_A}}$ in the WDVV equation $WDVV((i,1),(i,j-1),\mu_A,\mu_A)$ for $\F_{A}$
gives the same recursion relation with the one provided by 
the coefficient of the term $t^{\beta-e_{i,j}}e^{t_{\mu_{A_{i}}}}$ in $WDVV((i,1),(i,j-1),\mu_{A_{i}},\mu_{A_{i}})$ for $\F_{A_{i}}$
by the assumption and elementary calculation.
\qed
\end{pf}

\begin{lem}[Step 2]\label{atype decomp-2}
Assume that $c(\alpha,0)=c'(\alpha,0)$ and $c(\alpha'+\sum^{r}_{k\ne i} e_{k,1},1)=c'(\alpha',1)$
if $|\alpha|\le k+2$ and $|\alpha'|\le k$ for some $k\in \NN$.
Then we have $c(\beta +e_{i,j}+e_{i,j'}+e_{i,a_{i}-1},0)=c'(\beta +e_{i,j}+e_{i,j'}+e_{i,a_{i}-1},0)$ for some $i,j,j'$ if $|\beta|=k$.
\end{lem}
\begin{pf}
The coefficient of the term  $t^{\beta }(\prod_{k=2}^{r}t_{k,1})e^{t_{\mu_A }}$ in $WDVV((1,j),(1,j'),\mu_A,\mu_A)$ for $\F_{A}$ 
gives the same recursion relation with the one provided by the coefficient of the term $t^{\beta }e^{t_{\mu_A }}$
in the WDVV equation $WDVV((1,j),(1,j'),\mu_{A_{i}},\mu_{A_{i}})$ for $\F_{A_{i}}$ by the assumption and Lemma \ref{atype decomp-1}.
\qed
\end{pf}

\begin{lem}[Step 3]\label{atype decomp-3}
Assume that $c(\alpha,0)=c'(\alpha,0)$ and $c(\alpha'+\sum^{r}_{k\ne i} e_{k,1},1)=c'(\alpha',1)$
if $|\alpha|\le k+2$ and $|\alpha'|\le k$ for some $k\in \NN$.
Then we have $c(\gamma,0)=c'(\gamma,0)$ if $|\gamma|=k+3$ and $\gamma-e_{i,1}\ge 0$.
\end{lem}
\begin{pf}
We will show this claim by the induction on the degree of parameter $t_{i,j}$.
By Lemma \ref{atype decomp-2}, we have $c(\beta +e_{i,j}+e_{i,j'}+e_{i,a_{i}-1},0)=c'(\beta +e_{i,j}+e_{i,j'}+e_{i,a_{i}-1},0)$ 
if $|\beta|=k$.
Assume that $c(\gamma',0)=c'(\gamma',0)$ if $|\gamma'|=k+3$ and $\gamma'-e_{i,1}-e_{i,n}\geq 0$ for $n\ge l$.
Then we shall show that $c(\gamma,0)=c'(\gamma,0)$ for $|\gamma'|=k+3$ and $\gamma'-e_{i,1}-e_{i,l-1}\geq 0$.
We have $\deg(t^{\gamma-e_{i,1}-e_{i,l-1}})=l/a_{i}$. 
By the condition {\rm (iv)}, there exist $e_{i,m},e_{i,m'}$ such that 
\begin{itemize}
\item $\gamma-e_{i,1}-e_{i,l-1}-e_{i,m}-e_{i,m'}\geq 0$,
\item $\deg(t_{i,m})+\deg(t_{i,m'})\le l/a_{i}$.
\end{itemize}
Note that $\deg(t_{i,m})+\deg(t_{i,l})<1$.
We put $\beta:=\gamma-e_{i,l-1}-e_{i,m}-e_{i,m'}$.
Then the coefficient of the term $t^{\beta}(\prod_{k\ne i}t_{k,1})e^{t_{\mu_A}}$ in the WDVV equation $WDVV((i,m),(i,m'),l,\mu_A)$ for $\F_{A}$
gives the same recursion relation with the one provided by the coefficient of the term $t^{\beta}e^{t_{\mu_{A_{i}}}}$ in the WDVV equation 
$WDVV((i,m),(i,m'),l,\mu_{A_{i}})$ for $\F_{A_{i}}$ by the assumption and Lemma \ref{atype decomp-1}.
\qed
\end{pf}

\begin{lem}[Step 4]\label{atype decomp-4}
Assume that $c(\alpha,0)=c'(\alpha,0)$ and $c(\alpha'+\sum^{r}_{k\ne i} e_{k,1},1)=c'(\alpha',1)$
if $|\alpha|\le k+2$ and $|\alpha'|\le k$ for some $k\in \NN$.
Then we have $c(\gamma,0)=c'(\gamma,0)$ with $|\gamma|=k+3$.
\end{lem}
\begin{pf}
We will show this claim by the induction on the degree of parameter $t_{i,j}$.
By Lemma \ref{atype decomp-3}, we have $c(\gamma ,0)=c'(\gamma ,0)$ if $\gamma-e_{i,1}\ge 0$. 
Assume that $c(\gamma',0)=c'(\gamma',0)$ with $|\gamma'|=k+3$ and $\gamma'-e_{i,n}\geq 0$ for $n\leq l$.
We shall show that $c(\gamma,0)=c'(\gamma,0)$ with $|\gamma|\le k+3$ and $\gamma-e_{i,l+1}\ge 0$.
The coefficient of the term $t^{\gamma-e_{i,j}-e_{i,j'}-e_{i,l+1}}$ in the WDVV equation $WDVV((i,1),(i,l),(i,j),(i,j'))$ for $\F_{A}$
gives the same recursion relation with the one provided by the coefficient of the term
$t^{\gamma-e_{i,j}-e_{i,j'}-e_{i,l+1}}$ in the WDVV equation $WDVV((i,1),(i,l),(i,j),(i,j'))$ for $\F_{A_{i}}$
by the assumption.  
\qed
\end{pf}
Therefore we have Proposition~\ref{atype decomp-prop}.
\qed
\end{pf}

\subsection{Coefficients $c(\alpha,m)$ can be reconstructed}\label{reconst-s2}

In Subsection \ref{reconst-s1}, we showed that $c(\alpha,0)$ and $c(\alpha,1)$ can be reconstructed from $c(\beta,0)$ with $|\beta|=3$.
We define the total order $\prec$ on $\ZZ_{\ge 0}^{2}$ as follows:
\begin{itemize}
\item $(|\alpha|,m)\prec$ $(|\beta|,n)$ if $m<n$. 
\item $(|\alpha|,m)\prec$ $(|\beta|,m)$ if $|\alpha| <|\beta|$.
\end{itemize} 
We shall prove that $c(\alpha, m)$ can be reconstructed 
from $c(\beta,0)$ with $|\beta|=3$ by the induction on the well order $\prec$
on $\ZZ_{\ge 0}^{2}$.
\begin{propa}\label{mreconst}
A coefficient $c(\gamma,m)$ with $m\ge 2$ can be reconstructed from coefficients $c(\beta,0)$ with $|\beta|=3$.
\end{propa}
\begin{pf}
Assume that $c(\alpha, n)$ can be reconstructed from coefficients $c(\beta,0)$ with $|\beta|=3$ 
if $(|\alpha|,n)\prec (0,m-1)$.
First, we shall show that $c(0,m)$ can be reconstructed.  This coefficient must be zero for the case $\chi_A\le 0$.
For the case that $\chi_A>0$, we have the following Lemma~\ref{lem0-mreconst}:
\begin{lema}\label{lem0-mreconst}
A coefficient $c(0,m)$ can be reconstructed from coefficients $c(\alpha, n)$
with $(|\alpha|,n)\prec (0,m)$.
\end{lema}
\begin{pf}
We shall calculate the coefficient of $e^{mt_{\mu_A}}$
in $WDVV((i,1),(i,a_{i}-1),\mu_A,\mu_A)$. Then we have
\[
c(e_{i,1}+e_{i,a_{i}-1}+e_{1},0)\cdot 1 \cdot m^{3}\cdot c(0,m)
+ (known\ \ terms)=0.
\]
Therefore, the cofficient $c(0,m)$ can be reconstructed from coefficients $c(\alpha, n)$
satisfying $(|\alpha|,n)\prec (0,m)$.
\qed
\end{pf}
Next, we shall split the second step of the induction into following three cases. 
\begin{lema}[Case 1]\label{lem1-mreconst}
If a non--negative element $\gamma \in \ZZ^{\mu_{A}-2}$ satisfies that $|\gamma|=k+1$
and $\gamma-e_{i,j}\ge 0$ for some $j$ such that $j\ge 2$,
then the coefficient $c(\gamma,m)$ can be reconstructed from coefficients $c(\alpha, n)$
with $(|\alpha|,n)\prec (k+1,m)$.
\end{lema}
\begin{pf}
We shall calculate the coefficient of the term  
$t^{\gamma-e_{i,j}}e^{mt_{\mu_A}}$
in the WDVV equation $WDVV((i,1),(i,j-1),\mu_A,\mu_A)$. Then we have
\[
s_{1,j-1,a_{i}-j}\cdot c(e_{i,1}+e_{i,j-1}+e_{i,a_{i}-j},0)\cdot a_{i} \cdot m^{2}\cdot \gamma_{i,j} \cdot c(\gamma ,m)
+ (known\ \ terms)=0,
\]
Therefore, the coefficient $c(\gamma ,m)$ can be reconstructed from coefficients $c(\alpha, n)$
satisfying $(|\alpha|,n)\prec (k+1,m)$.
\qed
\end{pf}
\begin{lema}[Case 2]\label{lem2-mreconst}
If a non--negative element $\gamma \in \ZZ^{\mu_{A}-2}$ satisfies that $|\gamma|=k+1$
and $\gamma =\sum_{k=1}^{r}\gamma_{k,1}e_{k,1}$ for some  $\gamma_{1,1},\dots,\gamma_{r,1}$ such that $\prod_{k=1}^{r}\gamma_{k,1}\ne 0$,
then the coefficient $c(\gamma ,m)$ can be reconstructed from coefficients $c(\alpha, n)$
with $(|\alpha|,n)\prec (k+1,m)$.
\end{lema}
\begin{pf}
We shall calculate the coefficient of the term $(\prod_{k=1}^{r}t_{k,1}^{\gamma_{k,1}})e^{mt_{\mu_A}}$
in the WDVV equation $WDVV((i,1),(i,a_{i}-1),\mu_A,\mu_A)$. Then we have
\begin{eqnarray*}
\lefteqn{\hspace{-50mm}{\rm (i)} \ \{ c(e_{i,1}+e_{i,a_{i}-1}+e_{1},0)\cdot m^{3} +4\cdot 
c(2e_{i,1}+2e_{i,a_{i}-1},0)\cdot a_{i} 
\cdot m^{2}\cdot \gamma_{i,1}\} \cdot c(\gamma ,m)}\\ 
&&+ (known \ \ terms)=0 \\
&&\text{if} \ \ a_{i}\ge 3,\\
\lefteqn{\hspace{-50mm}{\rm (ii)} \ \{ 2c(2e_{i,1}+e_{1},0)\cdot m^{3} +24\cdot c(4e_{i,1},0)\cdot 2\cdot m^{2} 
\cdot \gamma_{i,1}\} \cdot c(\gamma ,m)}\\ 
&&+ (known \ \ terms)=0 \\
&&\text{if} \ \ a_{i}=2.
\end{eqnarray*}
If $\gamma_{i,1}\neq m$ for some $i$, the coefficient $c(\gamma,m)$ can be reconstructed from $c(\alpha, n)$
with $(\alpha,n)\prec (k+1,m)$.
If $\gamma_{i,1}=m$ for all $i$, we have $\deg ((\prod_{k=1}^{r}t_{k,1}^{\gamma_{k,1}})e^{mt_{\mu_A}})=2m\ge 2$ and hence $c(\gamma,m)=0$ except for the case $m=1$. 
\qed
\end{pf}
\begin{lema}[Case 3]\label{lem3-mreconst}
If a non--negative element $\gamma \in \ZZ^{\mu_{A}-2}$ satisfies that $|\gamma|=k+1$ and 
$\gamma =\sum_{k=1}^{r}\gamma_{k,1}e_{k,1}$ for some $\gamma_{1,1},\dots,\gamma_{r,1}$ such that $\prod_{k=1}^{r}\gamma_{k,1}= 0$,
then the coefficient $c(\gamma ,m)$ can be reconstructed from coefficients $c(\alpha, n)$
with $(|\alpha|,n)\prec (k+1,m)$.
\end{lema}
\begin{pf}
Assume that $\gamma_{i,1}=0$.
We shall calculate the coefficient of the term $(\prod_{k=1}^{r}t_{k,1}^{\gamma_{k,1}})e^{mt_{\mu_A}}$
in the WDVV equation $WDVV((i,1),(i,a_{i}-1),\mu_A,\mu_A)$. Then we have
\[
c(e_{1}+e_{i,a_{i}-1}+e_{i,1},0)\cdot m^{3} \cdot c(\gamma ,m)
+(known \ \ terms)=0.
\]
Therefore, the coefficient $c(\gamma ,m)$ can be reconstructed from coefficients $c(\alpha, n)$
satisfying $(|\alpha|,n)\prec (k+1,m)$.
\qed
\end{pf}
Hence, we have Proposition~\ref{mreconst}.
\qed
\end{pf}
We finish the proof of Theorem~\ref{first}.
\section{The condition \rm (iv)}

In the previous paper \cite{ist:1}, we showed Theorem \ref{first} for the case that $r=3$ and
derived the condition {\rm (iv)} by other conditions if $A$ satisfies $a_{2}\ge 3$. 
For the cases $A=(2,2,a_{3})$ with $a_{3}\ge 3$, we also derived the condition {\rm (iv)}
under the weaker condition than {\rm (iv)}:
\begin{itemize}
\item If $a_{i_1}=a_{i_2}$ for some $i_1,i_2\in\{1,2,3\}$, then the
Frobenius potential $\F$ is invariant under
the permutation of parameters $t_{i_1,j}$ and $t_{i_2,j}$ $(j=1,\dots, a_{i_1}-1)$.
\end{itemize}
For the case $A=(2,2,2)$, we have to assume the condition {(iv)} since, even under the weaker condition above, we can obtain a different Frobenius potential
which satisfies the conditions in Theorem \ref{first} except for the condition {\rm (iv)}.  

In the present paper, we can also obtain the similar result as in Proposition of \cite{ist:1}. 
\begin{defn}
We classify a multiplet of positive integers $A$ as follows and call it
\begin{enumerate}
\item a general multiplet if $a_1\ge 2$ and $a_2\ge 3$, 
\item a semi--general multiplet if $a_{1}=a_{2}=2$ and $a_{3}\ge 3$,
\item a non--general multiplet if $a_{1}=a_{2}=a_{3}=2$.
\end{enumerate}
\end{defn}
Under the above classification, the following Theorem \ref{sep} is the main result in this section.
All steps of the proof for Theorem \ref{sep} except for Sublemma \ref{s1,l2-l3.1} are common to the ones in our previous paper \cite{ist:1}. 
In order to make the proof self--contained, we shall include all details of arguments even if the arguments are common to the ones in \cite{ist:1}. 
\begin{thm*}[Theorem \ref{sep}]
Suppose that $A$ is a general multiplet.
For a non--negative $\beta\in \ZZ^{\mu_A-2}$, we have
\[
c\left(\beta+\sum_{k=1}^3e_{i_k,j_{k}},0\right)\ne 0
\]
only if $i_1=i_2=i_3$.
Suppose that $A$ is a semi--general multiplet.
For a non--negative $\beta\in \ZZ^{\mu_A-2}$, we have
\[
c\left(\beta+\sum_{k=1}^3e_{i_k,j_{k}},0\right)\ne 0
\]
only if $i_1=i_2=i_3$ under the following condition:
\begin{itemize}
\item[{\rm (iv')}] If $a_{i_1}=a_{i_2}$ for some $i_1,i_2\in\{1,\dots,r\}$, then the
Frobenius potential $\F$ is invariant under
the permutation of parameters $t_{i_1,j}$ and $t_{i_2,j}$ $(j=1,\dots, a_{i_1}-1)$.
\end{itemize}
\end{thm*}
\begin{rem}\label{(iv) to (iv')}
Before going to the proof of Theorem~\ref{sep}, we shall explain why the condition {\rm (iv)} implies the condition {\rm (iv')}.
We shall consider the permutation of index $(i_1,j)$ and $(i_{2},j)$ $(j=1,\dots, a_{i_1}-1)$ and denote it by $\rho$. 
Moreover, we shall take four non--negative elements 
$\alpha=\sum_{j}\alpha_{i_1,j}e_{i_1,j}$ and $\alpha'=\sum_{j}\alpha'_{i_2,j}e_{i_2,j}$
such that $\alpha_{i_1, j}=\alpha'_{i_2,j}$ for all $j$, and 
$\beta=\sum_{i,j}\beta_{i,j}e_{i,j}$, $\beta'=\sum_{i,j}\beta'_{i,j}e_{i,j}$
such that $\beta_{\rho (i,j)}=\beta'_{i,j}$ for all $i, j$. 
Then, in the proof of Theorem~\ref{first} under the condition {\rm (iv)},
one can sees, inductively, that the quadratic equation in front of $t^{\beta}e^{mt_{\mu_A}}$ in $WDVV(a,b,c,d)$ is exactly same 
with the one in front of $t^{\beta'}e^{mt_{\mu_A}}$ in $WDVV(\rho(a),\rho(b),\rho(c),\rho(d))$
since the coefficients of the trivial part for the Frobenius potential $\F$ are 
invariant under the permutations in the condition {\rm (iv')} (Proposition~\ref{lem3}).
Therefore one can have $c(\alpha, m)=c(\alpha', m)$ inductively in the proof of Theorem~\ref{first}
under the condition {\rm (iv)}.
\end{rem}
The following Lemma~\ref{wsep} is a part of the condition {\rm (iv)} in Theorem~\ref{first} and Theorem \ref{sep}.
However, for the later convenience, we shall show that Lemma~\ref{wsep} is derived from other conditions in Theorem~\ref{first}.
\begin{lema}\label{wsep}
Let $\gamma \in \ZZ^{\mu_{A}-2}$ be a non--negative element satisfying that 
$|\gamma|=4$ and $\gamma -e_{i_{1},j_{1}}-e_{i_{2},j_{2}}\ge 0$ for some $i_1, i_2$ such that $i_{1}\ne i_{2}$.
If $a_{i_1}\ge 3$, then we have $c(\gamma,0)=0$.
\end{lema}
\begin{pf}
Note that $c(\beta,0)=0$ if $|\beta|=3$ and $\beta-e_{i_{1},j_{1}}-e_{i_{2},j_{2}}\ge 0$ for some $i_1, i_2$ such that $i_{1}\ne i_{2}$ by Proposition \ref{lem3}.
We shall split the proof into the following two steps:

\vspace{5pt}
\noindent
Step 1: We shall consider the case that the term $t^{\gamma}$ has, as a factor,
$t_{i_1,j_1}$ for some $i_1$ and $j_1$ such that $a_{i_1}\ge 3$ and  $j_1\ge 2$.
We shall split, moreover, Step 1 into following four cases: 
\begin{enumerate}
\item The term $t^{\gamma}$ has $t_{i_1,j_1}, t_{i_2,j_2}, t_{i_3,j_3}$ as factors for some $j_1, j_2, j_3$ where $i_1, i_2, i_3$ are pairwise distinct.
\item The term $t^{\gamma}$ has $t_{i,j}$ and $t_{i,j'}$ as factors for each $i=i_{1}, i_{2}$ and some $j, j'$ where $i_1\ne i_2$.
\item The term $t^{\gamma}$ has $t_{i_1,j_1}, t_{i_1,j'_1}, t_{i_1,j''_1}$ and only $t_{i_2,j_2}$ as factors for some $i_1, i_2, j_1, j'_{1}, j''_{1}, j_{2}$ where $i_1 \ne i_2$.
\item The term $t^{\gamma}$ has, as factors, $t_{i_1,j_1}$ and $t_{i_{2},j}, t_{i_{2},j'}, t_{i_{2},j''}$ for some $i_1, i_2, j_1, j_2, j'_2, j''_2$ where $i_1 \ne i_2$.
\end{enumerate}
\vspace{5pt}
\noindent
Step 2: We shall consider the case that $t^{\gamma}$ does not have, as factors, $t_{i,j}$ for any $i$ and $j$ such that $j\ge 2$. 

\vspace{5pt}
\noindent
\underline{\it Step 1--{\rm (i)}. The term $t^{\gamma}$ has $t_{i_1,j_1}, t_{i_2,j_2}, t_{i_3,j_3}$ as factors for some $j_1, j_2, j_3$. }

\begin{sublema}[Step 1--(i)]\label{sl1-l1-wsep}
Let $\gamma \in \ZZ^{\mu_{A}-2}$ be a non--negative element such that $|\gamma|=4$ and 
$\gamma -e_{i_{1},j_{1}}-e_{i_{2},j_{2}}-e_{i_{3},j_{3}}\ge 0$ for distinct $i_{1}$, $i_{2}$, $i_{3}$ and some $j_1, j_2, j_3$.
If $a_{i_1}\ge 3$ and $j_{1}\ge 2$, then we have $c(\gamma,0)=0$. 
\end{sublema}
\begin{pf}
We shall calculate the coefficient of the term $t^{\gamma -e_{i_{1},j_{1}}-e_{i_{2},j_{2}}-e_{i_{3},j_{3}}}$ 
in the WDVV equation $WDVV((i_{3},j_{3}),(i_{2},j_{2}),(i_{1},j_{1}-1),(i_{1},1))$.
Then we have
\[
\gamma_{i_{1},j_{1}}\gamma_{i_{2},j_{2}}\gamma_{i_{3},j_{3}}\cdot c(\gamma,0)\cdot a_{i_{1}}\cdot 
s_{a_{i_{1}}-j_{1},j_{1}-1,1}\cdot c(e_{i_{1},a_{i_{1}}-j_{1}}+e_{i_{1},j_{1}-1}+e_{i_{1},1},0)=0.
\]
Hence we have $c(\gamma ,0)=0$.
\qed
\end{pf}

\vspace{5pt}
\noindent
\underline{\it Step 1--{\rm (ii)}.  The term
$t^{\gamma}$ has $t_{i,j}$ and $t_{i,j'}$ as factors for each $i=i_{1}, i_{2}$ and some $j, j'$.}
 
\begin{sublema}[Step 1--(ii)]\label{sl1-l2-wsep}
Let
$\gamma \in \ZZ^{\mu_A-2}$ be $e_{i_{1},j_{1}}+e_{i_{1},j'_{1}}+e_{i_{2},j_{2}}+e_{i_{2},j'_{2}}$ for some $i_1, i_2, j_1, j'_1, j_2, j'_2$ such that $i_{1}\neq i_{2}$. 
If $a_{i_{1}}\geq 3$ and $j'_{1}\ge 2$, then we have $c(\gamma,0)=0$.
\end{sublema}
\begin{pf}
We calculate the coefficient of $t_{i_2,j'_2}$ in $WDVV((i_{2},j_{2}),(i_{1},j_{1}),(i_{1},1),(i_{1},j'_{1}-1))$. 
Then we have
\begin{eqnarray*}
\lefteqn{{\rm (i)} \ \gamma_{i_{1},j_{1}}\gamma_{i_{2},j_{2}}\gamma_{i_{1},j'_{1}}\cdot c(\gamma,0)\cdot a_{i_{1}}\cdot 
s_{1,j'_{1}-1,a_{i_{1}}-j'_{1}}\cdot c(e_{i_{1},1}+e_{i_{1},j'_{1}-1}+e_{i_{1},a_{i_{1}}-j'_{1}},0)}\\
&&-(\gamma'_{i_{1},1}+1)(\gamma'_{i_{2},j_{2}}+1)(\gamma'_{i_{1},j_{1}+j'_{1}-1}+1)\cdot c(\gamma' +e_{i_{1},1}+e_{i_{2},j_{2}}+e_{i_{1},j_{1}+j'_{1}-1},0)
\cdot a_{i_{1}}\cdot \\
&&s_{j_1,j'_{1}-1,a_{i_{1}}+1-j_{1}-j'_{1}}\cdot c(e_{i_{1},j_{1}}+e_{i_{1},j'_{1}-1}+e_{i_{1},a_{i_{1}}+1-j_{1}-j'_{1}},0)=0\\ 
&&\text{if} \ \ 3\leq j_{1}+j'_{1}\leq a_{i_{1}} \ \text{and} \ \text{where} \ \gamma'=\gamma-e_{i_{1},j_{1}}-e_{i_{2},j_{2}}-e_{i_{1},j'_{1}}=e_{i_2,j'_2},\\ 
\lefteqn{{\rm (ii)} \ \gamma_{i_{1},j_{1}}\gamma_{i_{2},j_{2}}\gamma_{i_{1},j'_{1}}\cdot c(\gamma,0)\cdot a_{i_{1}}\cdot s_{1,j'_{1}-1,a_{i_{1}}-j'_{1}}\cdot
c(e_{i_{1},1}+e_{i_{1},j'_{1}-1}+e_{i_{1},a_{i_1}-j'_{1}},0)=0}\\
&&\text{if} \ \ j_{1}+j'_{1}> a_{i_{1}}.
\end{eqnarray*}
We shall show that $c(\gamma' +e_{i_{1},1}+e_{i_{2},j_{2}}+e_{i_{1},j_{1}+j'_{1}-1},0)=0$. We have the inequality $\deg(t_{i_{2},j_{2}}t_{i_{2},j'_{2}})=(j_{1}+j'_{1})/a_{i_{1}}\le 1$. 
We shall calculate the coefficient of the term $t_{i_{1},1}$ in the WDVV equation $WDVV((i_{2},j_{2}),(i_{2},j'_{2}),(i_{1},j_{1}+j'_{1}-2),(i_{1},1))$. 
Then we have
\begin{multline*}
(\gamma'_{i_{2},j_{2}}+1)(\gamma'_{i_{2},j'_{2}}+1)(\gamma'_{i_{1},j_{1}+j'_{1}-1}+1)\cdot c(\gamma' +e_{i_{1},1}+e_{i_{2},j_{2}}+e_{i_{1},j_{1}+j'_{1}-1},0)
\cdot a_{i_{1}} \cdot \\
s_{1,j_{1}+j'_{1}-2,a_{i_{1}}+1-j_{1}-j'_{1}}\cdot c(e_{i_{1},1}+e_{i_{1},j_{1}+j'_{1}-2}+e_{i_{1},a_{i_{1}}+1-j_{1}-j'_{1}},0)=0.
\end{multline*}
By Proposition \ref{lem3}, we have $c(\gamma' +e_{i_{1},1}+e_{i_{2},j_{2}}+e_{i_{1},j_{1}+j'_{1}-1},0)=0$
and hence $c(\gamma,0)=0$.
\qed
\end{pf}

\vspace{5pt}
\noindent
\underline{\it Step 1--{\rm (iii)}. The term
$t^{\gamma}$ has $t_{i_1,j_1}, t_{i_1,j'_1}, t_{i_1,j''_1}$ and only $t_{i_2,j_2}$ as factors for}

\noindent
\underline{some $i_1, i_2, j_1, j'_1, j''_1, j_2$ such that $i_1 \ne i_2$.}

\vspace{10pt}
We shall split Step 1--{\rm (iii)} into Case 1 ($j_2\ge 2$) and Case 2 ($j_2=1$).
\vspace{-5pt}
\begin{sublema}[Step 1--(iii)--Case 1]\label{sl2-l2-wsep}
Let $\gamma \in \ZZ^{\mu_{A}-2}$ be $e_{i_{1},j_{1}}+e_{i_{1},j'_{1}}+e_{i_{1},j''_{1}}+e_{i_{2},j_{2}}$ for some $i_1, i_2, j_1, j'_1, j''_1, j_2$ such that $i_{1}\neq i_{2}$. 
If $a_{i_{1}}\ge 3$, $a_{i_{2}}\ge 3$, $j'_{1}\ge 2$ and $j_{2}\ge 2$, then we have $c(\gamma,0)=0$.
\end{sublema}

\begin{pf}
We calculate the coefficient of $t_{i_1,j''_1}$ in $WDVV((i_{2},j_{2}),(i_{1},j_{1}),(i_{1},1),(i_{1},j'_{1}-1))$. 
Then we have
\begin{eqnarray*}
\lefteqn{{\rm (i)} \ \gamma_{i_{1},j_{1}}\gamma_{i_{2},j_{2}}\gamma_{i_{1},j'_{1}}\cdot c(\gamma,0)\cdot a_{i_{1}}\cdot 
s_{1,j'_{1}-1,a_{i_{1}}-j'_{1}}\cdot c(e_{i_{1},1}+e_{i_{1},j'_{1}-1}+e_{i_{1},a_{i_{1}}-j'_{1}},0)}\\
&&-(\gamma'_{i_{1},1}+1)(\gamma'_{i_{2},j_{2}}+1)(\gamma'_{i_{1},j_{1}+j'_{1}-1}+1)\cdot c(\gamma' +e_{i_{1},1}+e_{i_{2},j_{2}}+e_{i_{1},j_{1}+j'_{1}-1},0)
\cdot a_{i_{1}}\cdot \\
&&s_{j_1,j'_{1}-1,a_{i_{1}}+1-j_{1}-j'_{1}}\cdot c(e_{i_{1},j_{1}}+e_{i_{1},j'_{1}-1}+e_{i_{1},a_{i_{1}}+1-j_{1}-j'_{1}},0)=0\\ 
&&\text{if} \ \ 3\leq j_{1}+j'_{1}\leq a_{i_{1}} \ \text{and} \ \text{where} \ \gamma'=\gamma-e_{i_{1},j_{1}}-e_{i_{2},j_{2}}-e_{i_{1},j'_{1}}=e_{i_1,j''_1},\\
\lefteqn{{\rm (ii)} \ \gamma_{i_{1},j_{1}}\gamma_{i_{2},j_{2}}\gamma_{i_{1},j'_{1}}\cdot c(\gamma,0)\cdot a_{i_{1}}\cdot s_{1,j'_{1}-1,a_{i_{1}}-j'_{1}}\cdot
c(e_{i_{1},1}+e_{i_{1},j'_{1}-1}+e_{i_{1},a_{i_1}-j'_{1}},0)=0}\\
&&\text{if} \ \ j_{1}+j'_{1}> a_{i_{1}}.
\end{eqnarray*}
We shall show that $c(\gamma'+e_{i_{1},1}+e_{i_{2},j_{2}}+e_{i_{1},j_{1}+j'_{1}-1},0)=0$.
We have the inequality $\deg(t_{i_{1},j''_{1}})\leq (j_{1}+j'_{1}-1)/a_{i_{1}}$. We shall calculate the coefficient of the term 
$t_{i_{1},1}$ in the WDVV equation $WDVV((i_{1},j''_{1}),(i_{1},j_{1}+j'_{1}-1),(i_{2},j_{2}-1),(i_{2},1))$.
Then we have 
\begin{multline*}
\lefteqn{(\gamma'_{i_{1},j''_{1}}+1)(\gamma'_{i_{1},j_{1}+j'_{1}-1}+1)(\gamma'_{i_{2},j_{2}}+1)\cdot c(\gamma' +e_{i_{1},1}+e_{i_{2},j_{2}}+e_{i_{1},j_{1}+j'_{1}-1},0)
\cdot a_{i_{2}}\cdot}\\
s_{1,j_{2}-1,a_{i_{2}}-j_{2}}\cdot c(e_{i_{2},1}+e_{i_{2},j_{2}-1}+e_{i_{2},a_{i_{2}}-j_{2}},0)=0.
\end{multline*}
By Proposition \ref{lem3}, we have $c(\gamma'+e_{i_{1},1}+e_{i_{2},j_{2}}+e_{i_{1},j_{1}+j'_{1}-1},0)=0$
and hence $c(\gamma,0)=0$.
\qed
\end{pf}

\begin{sublema}[Step 1--(iii)--Case 2]\label{sl3-l2-wsep}
Let
$\gamma \in \ZZ^{\mu_{A}-2}$ be 
$e_{i_{1},j_{1}}+e_{i_{1},j'_{1}}+e_{i_{1},j''_{1}}+e_{i_{2},j_{2}}$ for  some $i_1, i_2, j_1, j'_1, j''_1, j_2$ such that $i_{1}\neq i_{2}$. 
If $a_{i_{1}}\ge 3$, $j'_{1}\ge 2$ and $j_{2}=1$, then we have $c(\gamma,0)=0$.
\end{sublema}
\begin{pf}
We calculate the coefficient of $t_{i_1,j''_1}$ in $WDVV((i_{2},j_{2}),(i_{1},j_{1}),(i_{1},1),(i_{1},j'_{1}-1))$. 
Then we have
\begin{eqnarray*}
\lefteqn{{\rm (i)} \ \gamma_{i_{1},j_{1}}\gamma_{i_{2},j_{2}}\gamma_{i_{1},j'_{1}}\cdot c(\gamma,0)\cdot a_{i_{1}}\cdot 
s_{1,j'_{1}-1,a_{i_{1}}-j'_{1}}\cdot c(e_{i_{1},1}+e_{i_{1},j'_{1}-1}+e_{i_{1},a_{i_{1}}-j'_{1}},0)}\\
&&-(\gamma'_{i_{1},1}+1)(\gamma'_{i_{2},j_{2}}+1)(\gamma'_{i_{1},j_{1}+j'_{1}-1}+1)\cdot c(\gamma' +e_{i_{1},1}+e_{i_{2},j_{2}}+e_{i_{1},j_{1}+j'_{1}-1},0)
\cdot a_{i_{1}}\cdot \\
&&s_{j_1,j'_{1}-1,a_{i_{1}}+1-j_{1}-j'_{1}}\cdot c(e_{i_{1},j_{1}}+e_{i_{1},j'_{1}-1}+e_{i_{1},a_{i_{1}}+1-j_{1}-j'_{1}},0)=0\\ 
&&\text{if} \ \ 3\leq j_{1}+j'_{1}\leq a_{i_{1}} \ \text{and} \ \text{where} \ \gamma'=\gamma-e_{i_{1},j_{1}}-e_{i_{2},j_{2}}-e_{i_{1},j'_{1}}=e_{i_1,j''_1},\\ 
\lefteqn{{\rm (ii)} \ \gamma_{i_{1},j_{1}}\gamma_{i_{2},j_{2}}\gamma_{i_{1},j'_{1}}\cdot c(\gamma,0)\cdot a_{i_{1}}\cdot s_{1,j'_{1}-1,a_{i_{1}}-j'_{1}}\cdot
c(e_{i_{1},1}+e_{i_{1},j'_{1}-1}+e_{i_{1},a_{i_1}-j'_{1}},0)=0}\\
&&\text{if} \ \ j_{1}+j'_{1}> a_{i_{1}}.
\end{eqnarray*}
We shall show that $c(\gamma' +e_{i_{1},1}+e_{i_{2},j_{2}}+e_{i_{1},j_{1}+j'_{1}-1},0)=0$.
We have the inequality $\deg(t_{i_{1},j''_{1}})\leq (j_{1}+j'_{1}-1)/a_{i_{1}}$. 
We shall calculate the coefficient of the term $t_{i_{1},1}$ in the WDVV equation 
$WDVV((i_{1},j''_{1}),(i_{2},j_{2}),(i_{1},j_{1}+j'_{1}-2),(i_{1},1))$. Then we have
\begin{eqnarray*}
\lefteqn{{\rm (i)} \ (\gamma'_{i_{1},j''_{1}}+1)(\gamma'_{i_{2},j_{2}}+1)(\gamma'_{i_{1},j_{1}+j'_{1}-1}+1)\cdot 
c(\gamma' +e_{i_{1},1}+e_{i_{2},j_{2}}+e_{i_{1},j_{1}+j'_{1}-1},0)\cdot a_{i_{1}}\cdot }\\
&&s_{1,j_{1}+j'_{1}-2,a_{i_{1}}+1-j_{1}-j'_{1}}\cdot c(e_{i_{1},1}+e_{i_{1},j_{1}+j'_{1}-2}+e_{i_{1},a_{i_{1}}+1-j_{1}-j'_{1}},0)\\
&&-\ \ (\gamma''_{i_{1},1}+1)(\gamma''_{i_{2},j_{2}}+1)(\gamma''_{i_{1},a_{i_{1}}-1}+1)\cdot 
c(\gamma' +2e_{i_{1},1}-e_{i_{1},j''_{1}}+e_{i_{2},j_{2}}+e_{i_{1},a_{i_{1}}-1},0)\cdot a_{i_{1}}\cdot \\
&&s_{j_{1}+j'_{1}-1,j''_{1},1}\cdot c(e_{i_{1},j_{1}+j'_{1}-2}+e_{i_{1},j''_{1}}+e_{i_{1},1},0)\\
&&\text{if} \ \ \deg(t_{i_{1},j''_{1}})=\frac{j_{1}+j'_{1}-1}{a_{i_{1}}} \ \text{and} \ \text{where} \ \gamma''=\gamma' +e_{i_{1},1}-e_{i_{1},j''_{1}}=e_{i_1,1},\\
\lefteqn{{\rm (ii)} \ (\gamma'_{i_{1},j''_{1}}+1)(\gamma'_{i_{2},j_{2}}+1)(\gamma'_{i_{1},j_{1}+j'_{1}-1}+1)\cdot 
c(\gamma' +e_{i_{1},1}+e_{i_{2},j_{2}}+e_{i_{1},j_{1}+j'_{1}-1},0)\cdot a_{i_{1}}\cdot }\\
&&s_{1,j_{1}+j'_{1}-2,a_{i_{1}}+1-j_{1}-j'_{1}}\cdot c(e_{i_{1},1}+e_{i_{1},j_{1}+j'_{1}-2}+e_{i_{1},a_{i_{1}}+1-j_{1}-j'_{1}},0)\\
&&\text{if} \ \ \deg(t_{i_{1},j''_{1}})\leq \frac{j_{1}+j'_{1}-2}{a_{i_{1}}}.
\end{eqnarray*}
If we have $c(\gamma' +2e_{i_{1},1}-e_{i_{1},j''_{1}}+e_{i_{2},j_{2}}+e_{i_{1},a_{i_{1}}-1},0 )\neq 0$, we should have 
\[
2\deg(t_{i_{1},1})+\deg(t_{i_{1},a_{i_{1}}-1})+\deg(i_{2},j_{2})\leq 2 \displaystyle\Leftrightarrow  \frac{1}{a_{i_{1}}}+\frac{1}{a_{i_{2}}} \geq 1.
\]
This inequality contradicts the assumption that $a_{i_{1}}\ge 3$ and $a_{i_{2}}\ge 2$.
Then we have $c(\gamma' +2e_{i_{1},1}-e_{i_{1},j''_{1}}+e_{i_{2},j_{2}}+e_{i_{1},a_{i_{1}}-1},0 )=0$ and 
$c(\gamma' +e_{i_{1},1}+e_{i_{2},j_{2}}+e_{i_{1},j_{1}+j'_{1}-1},0)=0$,
and hence $c(\gamma ,0)=0$.
\qed
\end{pf}

\vspace{5pt}
\noindent
\underline{\it Step 1--{\rm (iv)}. The tem
$t^{\gamma}$ has, as factors, $t_{i_1,j_1}$ and $t_{i_{2},j_2}, t_{i_{2},j'_2}, t_{i_{2},j''_2}$ for} 

\noindent
\underline{some $i_1, i_2, j_1, j_2, j'_{2}, j''_{2}$ such that $i_1 \ne i_2$.}

\vspace{5pt}
If $a_{i_2}\ge 3$ and some $j_2\ge 2$, these cases are already dealt with in previous arguments 
in Step 1--{\rm (iii)}. If $a_{i_2}\ge 3$ and $j_2=j'_2=j''_2=1$, we have $\deg(t^{\gamma})>2$.
Therefore, we only have to consider the case $a_{i_2}=2$:
\vspace{-5pt}
\begin{sublema}[Step 1--{\rm (iv)}]\label{sl4-l2-wsep}
Let $\gamma \in \ZZ^{\mu_{A}-2}$ be a non--negative element satisfying that $|\gamma|=4$,
$\gamma_{i_{3},j_{3}}=0$ for all $j_{3}\ge 1$ and
$\gamma-e_{i_{1},j_{1}}-2e_{i_{2},1}\ge 0$ for pairwise distinct $i_1, i_2, i_3$ and some $j_1$ such that $j_{1}\ge 2$.
If $a_{i_{1}}\ge 3$ and $a_{i_{2}}=2$, then we have $c(\gamma,0)=0$. 
\end{sublema}
\begin{pf}
We shall calculate the coefficient of the term $t^{\gamma-e_{i_{1},j_{1}}-2e_{i_{2},1}}$ 
in the WDVV equation $WDVV((i_{2},1),(i_{2},1),(i_{1},j_{1}-1),(i_{1},1))$. Then 
we have
\[
\gamma_{i_{1},j_{1}}\gamma_{i_{2},1}\cdot (\gamma_{i_{2},1}-1)c(\gamma, 0)\cdot a_{i_{1}}\cdot 
c(e_{i_{1},j_{1}-1}+e_{i_{1},1}+e_{i_{1},a_{i_{1}}-j_{1}},0)=0.
\]
Hence we have $c(\gamma,0)=0$. 
\qed
\end{pf}

\noindent
\vspace{5pt}
\underline{\it Step 2. The term 
$t^{\gamma}$ does not have, as factors, $t_{i,j}$ 
for any $i$ and $j$ such that $j\ge 2$.}

\vspace{5pt}
We shall split Step 2 into the following two cases:

\noindent
Case 1: The term $t^{\gamma}$ has $t_{i_1,1}, t_{i_2,1}, t_{i_3,1}$ as factors where $i_1,i_2,i_3$ are pairwise distinct.

\noindent
Case 2: The term $t^{\gamma}$ has, as factors, only two parameters $t_{i_1,1}, t_{i_2,1}$
for some $i_1, i_2$ such that $i_1\ne i_2$. 
\vspace{-5pt}
\begin{sublema}[Step 2--Case 1]\label{sl2-l1-wsep}
Let $\gamma \in \ZZ^{\mu_{A}-2}$ be a non--negative element satisfying that 
$|\gamma|=4$ and $\gamma=\gamma_{i_1,1}e_{i_1,1}+\gamma_{i_2,1}e_{i_2,1}+\gamma_{i_3,1}e_{i_3,1}+\gamma_{i_4,1}e_{i_4,1}$
for some $\gamma_{i_1,1}, \gamma_{i_2,1}, \gamma_{i_3,1}$ such that $\gamma_{i_1,1}\gamma_{i_2,1}\gamma_{i_3,1}\ne 0$.
If $a_{i_{1}}\ge 3$, then we have $c(\gamma,0)=0$. 
\end{sublema}

\begin{pf}
By the assumption that $a_{i_1}\ge 3$,
we have the inequality:
\[
\deg (t^{\gamma})\ge 4\frac{a_{l}-1}{a_{l}}\geq 2, 
\]
where $a_{l}=\min \{a_{i_1},a_{i_2},a_{i_3},a_{i_4}\}$. 
The first equality is attained if and only if $a_{i_1}=a_{i_2}=a_{i_3}=a_{i_4}$.
If $a_{i_1}=a_{i_2}=a_{i_3}=a_{i_4}$, one also has $\deg (t^{\gamma})> 2$.
Hence we have $c(\gamma ,0)=0$.
\qed
\end{pf}

\begin{sublema}[Step 2--Case 2]\label{sl5-l2-wsep}
Let
$\gamma \in \ZZ^{\mu_{A}-2}$ be a non--negative element satisfying that $|\gamma|=4$ and
$\gamma=\gamma_{i_{1},1}e_{i_{1},1}+\gamma_{i_{2},1}e_{i_{2},1}$ for some $i_1, i_2$ such that $i_1\ne i_2$ 
and  some  $\gamma_{i_{1},1}, \gamma_{i_{2},1}$ such that $\gamma_{i_{1},1}\gamma_{i_{2},1}\ne 0$.
If $a_{i_{1}}\ge 3$, then we have $c(\gamma,0)=0$. 
\end{sublema}
\begin{pf}
We have the inequality:
\[
\deg (t^{\gamma})\ge 4\frac{a_{l}-1}{a_{l}}\geq 2, 
\]
where $a_{l}=\min \{a_{i_{1}},a_{i_{2}}\}$.
The first equality is attained if and only if $a_{i_1}=a_{i_2}$.
We also have $\deg (t^{\gamma})> 2$ if $a_{i_1}=a_{i_2}$.
Hence we have $c(\gamma,0)=0$.   
\qed
\end{pf}
Therefore we have Lemma~\ref{wsep}.
\qed
\end{pf}

\begin{propd}\label{lem3.1}
If $A$ is a general multiplet,
a coefficient $c(\alpha,1)$ with $|\alpha|\le r$ 
is none-zero if and only if $\alpha=\sum_{k=1}^{r}e_{k,1}$.
If $A$ is a semi--general multiplet,
we have a coefficient $c(\alpha,1)$ with $|\alpha|\le r$ 
is nonezero if and only if $\alpha=\sum_{k=1}^{r}e_{k,1}$ under the condition {\rm (iv')}.
In particular, we have $c(\sum_{k=1}^{r}e_{k,1},1)=1$  
by the condition {\rm (vi)} of Theorem~\ref{first}.
\end{propd}
\begin{pf}

We shall split the proof into following two cases.

\begin{lema}[Case 1]\label{lem2-lem3.1}
Suppose that $A$ is a general multiplet.
If a non--negative element $\gamma \in \ZZ^{\mu_{A}-2}$ satisfies that
$|\gamma|=r$ and $\gamma =\sum_{k=1}^{r}\gamma_{k,1}e_{k,1}$ 
for some $\gamma_{1,1}, \dots, \gamma_{r,1}$ such that $\prod_{k=1}^{r}\gamma_{k,1}=0$,
then we have $c(\gamma ,1)=0$.
\end{lema}
\begin{pf}
Note that $c(\alpha,0)=0$ if $|\alpha|=4$ and $\alpha-e_{i_{1},j_{1}}-e_{i_{2},j_{2}}\ge 0$ for $i_{1}\neq i_{2}$
by Lemma~\ref{wsep} and $c(\alpha, 1)=0$ if $|\alpha|\le r-1$
since $\deg(t^{\alpha}e^{t_{\mu_{A}}})<2$. Assume that $\gamma_{i,1}=0$.
We shall calculate the coefficient of the term $(\prod_{k=1}^{r}t_{k,1}^{\gamma_{k,1}})e^{t_{\mu_A}}$
in $WDVV((i,1),(i,a_{i}-1),\mu_A,\mu_A)$. Then we have
\[
c(e_{1}+e_{i,a_{i}-1}+e_{i,1},0)\cdot c(\gamma ,1)=0.
\]
Hence we have $c(\gamma ,1)=0$ and then Lemma~\ref{lem2-lem3.1}.
\qed
\end{pf}
\begin{lemd}[Case 2]\label{s1,l2-l3.1}
Suppose that $A$ is a semi--general multiplet.
If a non--negative element $\gamma \in \ZZ^{\mu_{A}-2}$ satisfies that
$|\gamma|=r$ and $\gamma =\sum_{k=1}^{r}\gamma_{k,1}e_{k,1}$ 
for some $\gamma_{1,1}, \dots, \gamma_{r,1}$ such that $\prod_{k=1}^{r}\gamma_{k,1}=0$,
then we have $c(\gamma ,1)=0$.
\end{lemd}
\begin{pf}
Note that $c(\alpha,0)=0$ if $|\alpha|=4$, $\alpha-e_{i,1}-e_{i_{3},j_{3}}\ge 0$ for $1\le j_{3}\le a_{i_{3}}-1$, $i=1,2$
and $i_{3}\ge 3$ by Lemma \ref{wsep}, and that $c(\alpha, 1)=0$ if $|\alpha|\le 2$ since $\deg(t^{\alpha}e^{t_{\mu_{A}}})<2$.
If $\gamma_{i_{3}}=0$ for $i_{3}\ge3$, we have $c(\gamma,1)=0$ by the same argument in Lemma~\ref{lem2-lem3.1}.
Then it is enough to consider the following two cases:
\begin{enumerate}
\item $\gamma_{i_{3},1}=1$ for all $i_{3}\ge 3$, i.e., $\gamma=2e_{1,1}+\sum_{i=3}^{r}e_{i,1}$ or $\gamma=2e_{2,1}+\sum_{i=3}^{r}e_{i,1}$,
\item otherwise.
\end{enumerate}
For the case {\rm (ii)}, we have ${\rm deg}(t^{\gamma}e^{t_{\mu_A}})>2$ by easy argument.
Then we only have to consider the case {\rm (i)}. Without loss of generality, We assume that $\gamma=2e_{1,1}+\sum_{i=3}^{r}e_{i,1}$.
we shall calculate the coefficient of the term $(\prod_{i=3}^{r} t_{i,1})e^{t_{\mu_A}}$ in $WDVV((1,1),(1,1),(2,1),(2,1))$. Then we have
\[
2c(2e_{1,1}+e_{1},0)\cdot 1\cdot 2c(2e_{2,1}+\sum_{i=3}^{r}e_{i,1},1)+2c(2e_{2,1}+e_{1},0)\cdot 1\cdot 2c(2e_{1,1}+\sum_{i=3}^{r}e_{i,1},1)=0.
\]
we have $c(2e_{1,1}+\sum_{i=3}^{r}e_{i,1},1)=c(2e_{2,1}+\sum_{i=3}^{r}e_{i,1},1)$ by the condition {\rm (iv')} in Theorem~\ref{first}.
Hence we have $c(2e_{1,1}+\sum_{i=3}^{r}e_{i,1},1)=c(2e_{2,1}+\sum_{i=3}^{r}e_{i,1},1)=0$. 
\qed
\end{pf}
Therefore we have Proposition~\ref{lem3.1}.
\qed
\end{pf}
\begin{thm}\label{sep}
Suppose that $A$ is a general multiplet.
For a non--negative $\beta\in \ZZ^{\mu_A-2}$, we have
\[
c\left(\beta+\sum_{k=1}^3e_{i_k,j_{k}},0\right)\ne 0
\]
only if $i_1=i_2=i_3$.
Suppose that $A$ is a semi--general multiplet.
For a non--negative $\beta\in \ZZ^{\mu_A-2}$, we have
\[
c\left(\beta+\sum_{k=1}^3e_{i_k,j_{k}},0\right)\ne 0
\]
only if $i_1=i_2=i_3$ under the following condition:
\begin{itemize}
\item[{\rm (iv')}] If $a_{i_1}=a_{i_2}$ for some $i_1,i_2\in\{1,\dots,r\}$, then the
Frobenius potential $\F$ is invariant under
the permutation of parameters $t_{i_1,j}$ and $t_{i_2,j}$ $(j=1,\dots, a_{i_1}-1)$.
\end{itemize}
\end{thm}
\begin{pf}
We will prove Theorem~\ref{sep} by the induction on the length.
By Proposition \ref{lem3}, we have $c(\alpha,0)=0$ if $|\alpha|=3$ and $\alpha-e_{i_{1},j_{1}}-e_{i_{2},j_{2}}\geq 0$ for $i_{1}\neq i_{2}$.
Assume that $c(\alpha,0)=0$ if $|\alpha|\le k+3$ and $\alpha-e_{i_{1},j_{1}}-e_{i_{2},j_{2}}\geq 0$ for $i_{1}\neq i_{2}$.
Under this assumption,
we will prove that $c(\gamma,0)=0$ if $|\gamma|=k+4$ and
$\gamma-e_{i_{1},j_{1}}-e_{i_{2},j_{2}}\geq 0$ for $i_{1}\neq i_{2}$. 

We shall split the proof into the following four steps:

\vspace{5pt}
\noindent
Step 1: We shall consider the case that the term $t^{\gamma}$ has, as a factor,
$t_{i_1,j_1}$ for some $i_1$ and $j_1$ such that $a_{i_1}\ge 3$ and  $j_1\ge 2$.
We shall split, moreover, Step 1 into following four cases: 
\begin{enumerate}
\item The term $t^{\gamma}$ has $t_{i_1,j_1}, t_{i_2,j_2}, t_{i_3,j_3}$ as factors for some $j_1, j_2, j_3$ where $i_1, i_2, i_3$ are pairwise distinct.
\item The term $t^{\gamma}$ has $t_{i,j}$ and $t_{i,j'}$ as factors for each $i=i_{1}, i_{2}$ and some $j, j'$ where $i_1\ne i_2$.
\item The term $t^{\gamma}$ has $t_{i_1,j_1}, t_{i_1,j'_1}, t_{i_1,j''_1}$ and only $t_{i_2,j_2}$ as factors for some $i_1, i_2, j_1, j'_{1}, j''_{1}, j_{2}$ where $i_1 \ne i_2$.
\item The term $t^{\gamma}$ has, as factors, $t_{i_1,j_1}$ and $t_{i_{2},j}, t_{i_{2},j'}, t_{i_{2},j''}$ for some $i_1, i_2, j_1, j_2, j'_2, j''_2$ where $i_1 \ne i_2$.
\end{enumerate}

\vspace{5pt}
\noindent
Step 2: We shall consider the case that $a_{i_1}\ge 3$ and $t^{\gamma}$ does not have, as factors, $t_{i,j}$ for 
any $i$ and $j$ such that $j\ge 2$. 

\vspace{5pt}
\noindent
Step 3: We shall consider the case that $a_{1}=a_{2}=2$ and $a_3\ge 3$, i.e., $A$ is a semi--general multiplet.

\vspace{5pt}
\noindent
\underline{\it Step 1--{\rm (i)}. The term $t^{\gamma}$ has $t_{i_1,j_1}, t_{i_2,j_2}, t_{i_3,j_3}$ as factors for some $j_1, j_2, j_3$. }

\begin{sublema}[Step 1--(i)]\label{sl1-l1-sep}
Assume that $a_{i_{1}}\ge 3$ and $c(\alpha,0)=0$ 
if $|\alpha|\le k+3$ and $\alpha-e_{i,j}-e_{i',j'}\geq 0$ for some  $i,i',j,j'$ such that $i\neq i'$.
If a non--negative element $\gamma \in \ZZ^{\mu_{A}-2}$ satisfies that 
$|\gamma|=k+4$, $\gamma -e_{i_{1},j_{1}}-e_{i_{2},j_{2}}-e_{i_{3},j_{3}}\ge 0$ for pairwise distinct $i_{1}$, $i_{2}$, $i_{3}$ and some $j_1$ such that $j_{1}\ge 2$,
then we have $c(\gamma,0)=0$. 
\end{sublema}
\begin{pf}
We shall calculate the coefficient of the term $t^{\gamma -e_{i_{1},j_{1}}-e_{i_{2},j_{2}}-e_{i_{3},j_{3}}}$ in the WDVV equation 
$WDVV((i_{3},j_{3}),(i_{2},j_{2}),(i_{1},j_{1}-1),(i_{1},1))$.
Then we have
\[
\gamma_{i_{1},j_{1}}\gamma_{i_{2},j_{2}}\gamma_{i_{3},j_{3}}\cdot c(\gamma,0)\cdot a_{i_{1}}\cdot 
s_{a_{i_{1}}-j_{1},j_{1}-1,1}\cdot c(e_{i_{1},a_{i_{1}}-j_{1}}+e_{i_{1},j_{1}-1}+e_{i_{1},1},0)=0.
\]
Hence we have $c(\gamma ,0)=0$.
\qed
\end{pf}

\vspace{5pt}
\noindent
\underline{\it Step 1--{\rm (ii)}.  The term
$t^{\gamma}$ has $t_{i,j}$ and $t_{i,j'}$ as factors for each $i=i_{1}, i_{2}$ and some $j, j'$.}

\begin{sublema}[Step 1--(ii)]\label{sl1-l2-sep}
Assume that $a_{i_{1}}\ge 3$ and $c(\alpha,0)=0$ 
if $|\alpha|\le k+3$ and $\alpha-e_{i,j}-e_{i',j'}\geq 0$ for some  $i,i',j,j'$ such that $i\neq i'$. 
If a non--negative element $\gamma \in \ZZ^{\mu_{A}-2}$ satisfies that $|\gamma|=k+4$,
$\gamma_{i_{3},j_{3}}=0$ for all $i_3\ne i_1,i_2$ and all $j_{3}\geq 1$ and that
$\gamma-e_{i_{1},j_{1}}-e_{i_{1},j'_{1}}-e_{i_{2},j_{2}}-e_{i_{2},j'_{2}}\geq 0$ for some $i_1,i_2,j_1,j'_1,j_2,j'_2$ such that $i_{1}\neq i_{2}$ and $j'_{1}\ge 2$, 
then we have $c(\gamma,0)=0$.
\end{sublema}
\begin{pf}
We shall calculate the coefficient of the term $t^{\gamma-e_{i_{1},j_{1}}-e_{i_{2},j_{2}}-e_{i_{1},j'_{1}}}$ in the WDVV equation
$WDVV((i_{2},j_{2}),(i_{1},j_{1}),(i_{1},1),(i_{1},j'_{1}-1))$. 
Then we have
\begin{eqnarray*}
\lefteqn{{\rm (i)} \ \gamma_{i_{1},j_{1}}\gamma_{i_{2},j_{2}}\gamma_{i_{1},j'_{1}}\cdot c(\gamma,0)\cdot a_{i_{1}}\cdot 
s_{1,j'_{1}-1,a_{i_{1}}-j'_{1}}\cdot c(e_{i_{1},1}+e_{i_{1},j'_{1}-1}+e_{i_{1},a_{i_{1}}-j'_{1}},0)}\\
&&-(\gamma'_{i_{1},1}+1)(\gamma'_{i_{2},j_{2}}+1)(\gamma'_{i_{1},j_{1}+j'_{1}-1}+1)\cdot c(\gamma' +e_{i_{1},1}+e_{i_{2},j_{2}}+e_{i_{1},j_{1}+j'_{1}-1},0)
\cdot a_{i_{1}}\cdot \\
&&s_{j_1,j'_{1}-1,a_{i_{1}}+1-j_{1}-j'_{1}}\cdot c(e_{i_{1},j_{1}}+e_{i_{1},j'_{1}-1}+e_{i_{1},a_{i_{1}}+1-j_{1}-j'_{1}},0)=0\\ 
&&\text{if} \ \ 3\leq j_{1}+j'_{1}\leq a_{i_{1}} \ \text{and} \ \text{where} \ \gamma'=\gamma-e_{i_{1},j_{1}}-e_{i_{2},j_{2}}-e_{i_{1},j'_{1}},\\ 
\lefteqn{{\rm (ii)} \ \gamma_{i_{1},j_{1}}\gamma_{i_{2},j_{2}}\gamma_{i_{1},j'_{1}}\cdot c(\gamma,0)\cdot a_{i_{1}}\cdot s_{1,j'_{1}-1,a_{i_{1}}-j'_{1}}\cdot
c(e_{i_{1},1}+e_{i_{1},j'_{1}-1}+e_{i_{1},a_{i_1}-j'_{1}},0)=0}\\
&&\text{if} \ \ j_{1}+j'_{1}> a_{i_{1}}.
\end{eqnarray*}
We shall show that $c(\gamma' +e_{i_{1},1}+e_{i_{2},j_{2}}+e_{i_{1},j_{1}+j'_{1}-1},0)=0$.
we have the inequality $\deg(t^{\gamma'}t_{i_{2},j_{2}})=(j_{1}+j'_{1})/a_{i_{1}}$, i.e, $\deg(t_{i_{2},j_{2}}t_{i_{2},j'_{2}})\leq (j_{1}+j'_{1})/a_{i_{1}}\le 1$. 
We shall calculate the coefficient of the term $t^{\gamma' +e_{i_{1},1}-e_{i_{2},j'_{2}}}$ in $WDVV((i_{2},j_{2}),(i_{2},j'_{2}),(i_{1},j_{1}+j'_{1}-2),(i_{1},1))$. 
Then we have
\begin{multline*}
(\gamma'_{i_{2},j_{2}}+1)(\gamma'_{i_{2},j'_{2}}+1)(\gamma'_{i_{1},j_{1}+j'_{1}-1}+1)\cdot c(\gamma' +e_{i_{1},1}+e_{i_{2},j_{2}}+e_{i_{1},j_{1}+j'_{1}-1},0)
\cdot a_{i_{1}} \cdot \\
s_{1,j_{1}+j'_{1}-2,a_{i_{1}}+1-j_{1}-j'_{1}}\cdot c(e_{i_{1},1}+e_{i_{1},j_{1}+j'_{1}-2}+e_{i_{1},a_{i_{1}}+1-j_{1}-j'_{1}},0)=0.
\end{multline*}
By Proposition \ref{lem3}, $c(\gamma' +e_{i_{1},1}+e_{i_{2},j_{2}}+e_{i_{1},j_{1}+j'_{1}-1},0)=0$.
Hence we have $c(\gamma,0)=0$.
\qed
\end{pf}

\vspace{5pt}
\noindent
\underline{\it Step 1--{\rm (iii)}. The term
$t^{\gamma}$ has $t_{i_1,j_1}, t_{i_1,j'_1}, t_{i_1,j''_1}$ and $t_{i_2,j_2}$ as factors for}

\noindent
\underline{some $i_1, i_2, j_1, j'_1, j''_1, j_2$ such that $i_{1}\neq i_{2}$.}

\vspace{10pt}
We shall split Step 1--{\rm (iii)} into Case 1 ($j_2\ge 2$) and Case 2 ($j_2=1$).
\vspace{-5pt}
\begin{sublema}[Step 1--(iii)--Case 1]\label{sl2-l2-sep}
Assume that $a_{i_{1}}\ge 3$, $a_{i_{2}}\ge 3$  and $c(\alpha,0)=0$ 
if $|\alpha|\le k+3$ and $\alpha-e_{i,j}-e_{i',j'}\geq 0$ for some  $i,i',j,j'$ such that $i\neq i'$. 
If a non--negative element $\gamma \in \ZZ^{\mu_{A}-2}$ satisfies that $|\gamma|=k+4$,
$\gamma_{i_{3},j_{3}}=0$ for all $i_3\ne i_1,i_2$ and all $j_{3}\geq 1$ and that
$\gamma-e_{i_{1},j_{1}}-e_{i_{1},j'_{1}}-e_{i_{1},j''_{1}}-e_{i_{2},j_{2}}\geq 0$ for some $i_1,i_2,j_1,j'_1,j''_1,j_2$ such that $j'_{1}\ge 2$ and $j_{2}\ge 2$, 
then we have $c(\gamma,0)=0$.
\end{sublema}
\begin{pf}
We shall calculate the coefficient of the term $t^{\gamma-e_{i_{1},j_{1}}-e_{i_{2},j_{2}}-e_{i_{1},j'_{1}}}$ in
the WDVV equation $WDVV((i_{2},j_{2}),(i_{1},j_{1}),(i_{1},1),(i_{1},j'_{1}-1))$. 
Then we have
\begin{eqnarray*}
\lefteqn{{\rm (i)} \ \gamma_{i_{1},j_{1}}\gamma_{i_{2},j_{2}}\gamma_{i_{1},j'_{1}}\cdot c(\gamma,0)\cdot a_{i_{1}}\cdot 
s_{1,j'_{1}-1,a_{i_{1}}-j'_{1}}\cdot c(e_{i_{1},1}+e_{i_{1},j'_{1}-1}+e_{i_{1},a_{i_{1}}-j'_{1}},0)}\\
&&-(\gamma'_{i_{1},1}+1)(\gamma'_{i_{2},j_{2}}+1)(\gamma'_{i_{1},j_{1}+j'_{1}-1}+1)\cdot c(\gamma' +e_{i_{1},1}+e_{i_{2},j_{2}}+e_{i_{1},j_{1}+j'_{1}-1},0)
\cdot a_{i_{1}}\cdot \\
&&s_{j_1,j'_{1}-1,a_{i_{1}}+1-j_{1}-j'_{1}}\cdot c(e_{i_{1},j_{1}}+e_{i_{1},j'_{1}-1}+e_{i_{1},a_{i_{1}}+1-j_{1}-j'_{1}},0)=0\\ 
&&\text{if} \ \ 3\leq j_{1}+j'_{1}\leq a_{i_{1}} \ \text{and} \ \text{where} \ \gamma'=\gamma-e_{i_{1},j_{1}}-e_{i_{2},j_{2}}-e_{i_{1},j'_{1}},\\ 
\lefteqn{{\rm (ii)} \ \gamma_{i_{1},j_{1}}\gamma_{i_{2},j_{2}}\gamma_{i_{1},j'_{1}}\cdot c(\gamma,0)\cdot a_{i_{1}}\cdot s_{1,j'_{1}-1,a_{i_{1}}-j'_{1}}\cdot
c(e_{i_{1},1}+e_{i_{1},j'_{1}-1}+e_{i_{1},a_{i_1}-j'_{1}},0)=0}\\
&&\text{if} \ \ j_{1}+j'_{1}> a_{i_{1}}.
\end{eqnarray*}
We shall show that $c(\gamma' +e_{i_{1},1}+e_{i_{2},j_{2}}+e_{i_{1},j_{1}+j'_{1}-1},0)=0$.
We have the inequality $\deg(t_{i_{1},j''_{1}})\leq (j_{1}+j'_{1}-1)/a_{i_{1}}$. We shall calculate the coefficient of the term 
$t^{\gamma' +e_{i_{1},1}-e_{i_{1},j''_{1}}}$ in the WDVV equation $WDVV((i_{1},j''_{1}),(i_{1},j_{1}+j'_{1}-1),(i_{2},j_{2}-1),(i_{2},1))$.
Then we have 
\begin{multline*}
\lefteqn{(\gamma'_{i_{1},j''_{1}}+1)(\gamma'_{i_{1},j_{1}+j'_{1}-1}+1)(\gamma'_{i_{2},j_{2}}+1)\cdot c(\gamma 
+e_{i_{1},1}+e_{i_{2},j_{2}}+e_{i_{1},j_{1}+j'_{1}-1},0)
\cdot a_{i_{2}}\cdot}\\
s_{1,j_{2}-1,a_{i_{2}}-j_{2}}\cdot c(e_{i_{2},1}+e_{i_{2},j_{2}-1}+e_{i_{2},a_{i_{2}}-j_{2}},0)=0.
\end{multline*}
By Proposition \ref{lem3}, we have $c(\gamma' +e_{i_{1},1}+e_{i_{2},j_{2}}+e_{i_{1},j_{1}+j'_{1}-1},0)=0$ and hence $c(\gamma,0)=0$.
\qed
\end{pf}

\begin{sublema}[Step 1--(iii)--Case 2]\label{sl3-l2-sep}
Assume that $a_{i_{1}}\ge 3$ and $c(\alpha,0)=0$ 
if $|\alpha|\le k+3$ and $\alpha-e_{i,j}-e_{i',j'}\geq 0$ for some  $i,i',j,j'$ such that $i\neq i'$. 
If a non--negative element $\gamma \in \ZZ^{\mu_{A}-2}$ satisfies that $|\gamma|=k+4$,
$\gamma_{i_{3},j_{3}}=0$ for all $i_3\ne i_1,i_2$ and all $j_{3}\geq 1$ and that
$\gamma-e_{i_{1},j_{1}}-e_{i_{1},j'_{1}}-e_{i_{1},j''_{1}}-e_{i_{2},j_{2}}\geq 0$ for  some $i_1,i_2,j_1,j'_1,j''_1,j_2$ such that \noindent $j'_{1}\ge 2$ and $j_{2}=1$, 
then we have $c(\gamma,0)=0$.
\end{sublema}
\begin{pf}
We shall calculate the coefficient of the term $t^{\gamma-e_{i_{1},j_{1}}-e_{i_{2},j_{2}}-e_{i_{1},j'_{1}}}$ in the WDVV equation 
$WDVV((i_{2},j_{2}),(i_{1},j_{1}),(i_{1},1),(i_{1},j'_{1}-1))$. 
Then we have
\begin{eqnarray*}
\lefteqn{{\rm (i)} \ \gamma_{i_{1},j_{1}}\gamma_{i_{2},j_{2}}\gamma_{i_{1},j'_{1}}\cdot c(\gamma,0)\cdot a_{i_{1}}\cdot 
s_{1,j'_{1}-1,a_{i_{1}}-j'_{1}}\cdot c(e_{i_{1},1}+e_{i_{1},j'_{1}-1}+e_{i_{1},a_{i_{1}}-j'_{1}},0)}\\
&&-(\gamma'_{i_{1},1}+1)(\gamma'_{i_{2},j_{2}}+1)(\gamma'_{i_{1},j_{1}+j'_{1}-1}+1)\cdot c(\gamma' +e_{i_{1},1}+e_{i_{2},j_{2}}+e_{i_{1},j_{1}+j'_{1}-1},0)
\cdot a_{i_{1}}\cdot \\
&&s_{j_1,j'_{1}-1,a_{i_{1}}+1-j_{1}-j'_{1}}\cdot c(e_{i_{1},j_{1}}+e_{i_{1},j'_{1}-1}+e_{i_{1},a_{i_{1}}+1-j_{1}-j'_{1}},0)=0\\ 
&&\text{if} \ \ 3\leq j_{1}+j'_{1}\leq a_{i_{1}} \ \text{and} \ \text{where} \ \gamma'=\gamma-e_{i_{1},j_{1}}-e_{i_{2},j_{2}}-e_{i_{1},j'_{1}},\\ 
\lefteqn{{\rm (ii)} \ \gamma_{i_{1},j_{1}}\gamma_{i_{2},j_{2}}\gamma_{i_{1},j'_{1}}\cdot c(\gamma,0)\cdot a_{i_{1}}\cdot s_{1,j'_{1}-1,a_{i_{1}}-j'_{1}}\cdot
c(e_{i_{1},1}+e_{i_{1},j'_{1}-1}+e_{i_{1},a_{i_{1}}-j'_{1}},0)=0}\\
&&\text{if} \ \ j_{1}+j'_{1}> a_{i_{1}}.
\end{eqnarray*}
We shall show that $c(\gamma' +e_{i_{1},1}+e_{i_{2},j_{2}}+e_{i_{1},j_{1}+j'_{1}-1},0)=0$.
We have the inequality $\deg(t_{i_{1},j''_{1}})\leq (j_{1}+j'_{1}-1)/a_{i_{1}}$. We shall calculate the coefficient of the term $t^{\gamma' +e_{i_{1},1}-e_{i_{1},j''_{1}}}$
in the WDVV equation
$WDVV((i_{1},j''_{1}),(i_{2},j_{2}),(i_{1},j_{1}+j'_{1}-2),(i_{1},1))$. Then we have
\begin{eqnarray*}
\lefteqn{{\rm (i)} \ (\gamma'_{i_{1},j''_{1}}+1)(\gamma'_{i_{2},j_{2}}+1)(\gamma'_{i_{1},j_{1}+j'_{1}-1}+1)\cdot 
c(\gamma' +e_{i_{1},1}+e_{i_{2},j_{2}}+e_{i_{1},j_{1}+j'_{1}-1},0)\cdot a_{i_{1}}\cdot }\\
&&s_{1,j_{1}+j'_{1}-2,a_{i_{1}}+1-j_{1}-j'_{1}}\cdot c(e_{i_{1},1}+e_{i_{1},j_{1}+j'_{1}-2}+e_{i_{1},a_{i_{1}}+1-j_{1}-j'_{1}},0)\\
&&-\ \ (\gamma''_{i_{1},1}+1)(\gamma''_{i_{2},j_{2}}+1)(\gamma''_{i_{1},a_{i_{1}}-1}+1)\cdot 
c(\gamma' +2e_{i_{1},1}-e_{i_{1},j''_{1}}+e_{i_{2},j_{2}}+e_{i_{1},a_{i_{1}}-1},0)\cdot a_{i_{1}}\cdot \\
&&s_{j_{1}+j'_{1}-1,j''_{1},1}\cdot c(e_{i_{1},j_{1}+j'_{1}-2}+e_{i_{1},j''_{1}}+e_{i_{1},1},0)=0\\
&&\text{if} \ \ \deg(t_{i_{1},j''_{1}})=\frac{j_{1}+j'_{1}-1}{a_{i_{1}}} \ \text{and} \ \text{where} \ \gamma''=\gamma' +e_{i_{1},1}-e_{i_{1},j''_{1}},
\end{eqnarray*}
\begin{eqnarray*}
\lefteqn{{\rm (ii)} \ (\gamma'_{i_{1},j''_{1}}+1)(\gamma'_{i_{2},j_{2}}+1)(\gamma'_{i_{1},j_{1}+j'_{1}-1}+1)\cdot 
c(\gamma' +e_{i_{1},1}+e_{i_{2},j_{2}}+e_{i_{1},j_{1}+j'_{1}-1},0)\cdot a_{i_{1}}\cdot }\\
&&s_{1,j_{1}+j'_{1}-2,a_{i_{1}}+1-j_{1}-j'_{1}}\cdot c(e_{i_{1},1}+e_{i_{1},j_{1}+j'_{1}-2}+e_{i_{1},a_{i_{1}}+1-j_{1}-j'_{1}},0)=0\\
&&\text{if} \ \ \deg(t_{i_{1},j''_{1}})\leq \frac{j_{1}+j'_{1}-2}{a_{i_{1}}}.
\end{eqnarray*}
If $c(\gamma' +2e_{i_{1},1}-e_{i_{1},j''_{1}}+e_{i_{2},j_{2}}+e_{i_{1},a_{i_{1}}-1},0 )\neq 0$, 
we should have 
\[
2\deg(t_{i_{1},1})+\deg(t_{i_{1},a_{i_{1}}-1})+\deg(i_{2},j_{2})\leq 2 \Leftrightarrow  \frac{1}{a_{i_{1}}}+\frac{1}{a_{i_{2}}} \geq 1.
\]
This inequality contradicts the assumption that $a_{i_{1}}\ge 3$ and $a_{i_{2}}=2$.
Then we have $c(\gamma' +2e_{i_{1},1}-e_{i_{1},j''_{1}}+e_{i_{2},j_{2}}+e_{i_{1},a_{i_{1}}-1},0 )=0$ and 
$c(\gamma' +e_{i_{1},1}+e_{i_{2},j_{2}}+e_{i_{1},j_{1}+j'_{1}-1},0)=0$.
Hence we have $c(\gamma ,0)=0$.
\qed
\end{pf}

\vspace{5pt}
\noindent
\underline{\it Step 1--{\rm (iv)}. The tem
$t^{\gamma}$ has, as factors, $t_{i_1,j_1}$ and $t_{i_{2},j_2}, t_{i_{2},j'_2}, t_{i_{2},j''_2}$ for}

\noindent
\underline{some $i_1, i_2, j_1, j_2, j'_2, j''_2$ such that $i_{1}\neq i_{2}$.}

\vspace{5pt}
If $a_{i_2}\ge 3$ and some $j\ge 2$, these cases are already dealt with in previous arguments 
in Step 1--{\rm (iii)}. If $a_{i_2}\ge 3$ and $j_2=j'_2=j''_2=1$, we have $\deg(t^{\gamma})>2$.
Therefore, we only have to consider the case $a_{i_2}=2$:
\vspace{-5pt}
\begin{sublema}[Step 1--(iv)]\label{sl4-l2-sep}
Assume that $a_{i_{1}}\ge 3$, $a_{i_{2}}=2$ and $c(\alpha,0)=0$ 
if $|\alpha|\le k+3$ and $\alpha-e_{i_{1},j_{1}}-e_{i_{2},j_{2}}\geq 0$ for  some $i_1,i_2,j_1,j_2$ such that $i_{1}\neq i_{2}$. 
If a non--negative element $\gamma \in \ZZ^{\mu_{A}-2}$ satisfies that $|\gamma|=k+4$,
$\gamma_{i_{3},j_{3}}=0$ for all $i_3\ne i_1,i_2$ and all $j_{3}\geq 1$ and that
$\gamma-e_{i_{1},j_{1}}-2e_{i_{2},1}\ge 0$ for some $i_1,i_2,j_1$ such that $j_{1}\ge 2$,
then we have $c(\gamma,0)=0$. 
\end{sublema}
\begin{pf}
We shall calculate the coefficient of the term $t^{\gamma-e_{i_{1},j_{1}}-2e_{i_{2},1}}$ 
in the WDVV equation $WDVV((i_{2},1),(i_{2},1),(i_{1},j_{1}-1),(i_{1},1))$. 
Then we have
\[
\gamma_{i_{1},j_{1}}\gamma_{i_{2},1}\cdot (\gamma_{i_{2},1}-1)c(\gamma, 0)\cdot a_{i_{1}}\cdot 
c(e_{i_{1},j_{1}-1}+e_{i_{1},1}+e_{i_{1},a_{i_{1}}-j_{1}},0)=0.
\]
Hence we have $c(\gamma,0)=0$. 
\qed
\end{pf}

\vspace{5pt}
\noindent
\underline{\it Step 2. $a_{i_1}\ge 3$ and $t^{\gamma}$ does not have, as factors, $t_{i,j}$ 
for any $i$ and $j$ such that $j\ge 2$.}

\vspace{10pt}
We shall split Step 2 into the following two cases:

\noindent
Case 1: The term $t^{\gamma}$ has $t_{i_1,1}, t_{i_2,1}, t_{i_3,1}$ as factors where $i_1,i_2,i_3$ are pairwise distinct.

\noindent
Case 2: The term $t^{\gamma}$ has, as factors, only two parameters $t_{i_1,1}, t_{i_2,1}$
for some $i_1,i_2$ such that $i_1\ne i_2$. 
\vspace{-5pt}
\begin{sublema}[Step 2--Case 1]\label{sl2-l1-sep}
Assume that $a_{i_{1}}\ge 3$ and $c(\alpha,0)=0$ 
if $|\alpha|\le k+3$ and $\alpha-e_{i,j}-e_{i',j'}\geq 0$ for some  $i,i',j,j'$ such that $i\neq i'$. 
If a non--negative element $\gamma \in \ZZ^{\mu_{A}-2}$ satisfies that 
$|\gamma|=k+4$ and 
$\gamma=\sum_{i=1}^{r}\gamma_{i,1}e_{i,1}$ such that  $\gamma_{i_1,1}\gamma_{i_2,1}\gamma_{i_3,1}\ne 0$ where $i_1,i_2,i_3$ are pairwise distinct,
then we have $c(\gamma,0)=0$. 
\end{sublema}
\begin{pf}
We have 
\[
\deg (t^{\gamma})\ge (k+4)\frac{a_{l}-1}{a_{l}}>2, 
\]
where $a_{l}=\min \{a_{1},\dots,a_{r}\}$. 
Hence we have $c(\gamma ,0)=0$.
\qed
\end{pf}
\begin{sublema}[Step 2--Case 2]\label{sl5-l2-sep}
Assume that $a_{i_{1}}\ge 3$ and $c(\alpha,0)=0$ 
if $|\alpha|\le k+3$ and $\alpha-e_{i,j}-e_{i',j'}\geq 0$ for some  $i,i',j,j'$ such that $i\neq i'$. 
If a non--negative element $\gamma \in \ZZ^{\mu_{A}-2}$ satisfies that $|\gamma|=k+4$
and $\gamma=\gamma_{i_{1},1}e_{i_{1},1}+\gamma_{i_{2},1}e_{i_{2},1}$ such that $i_1\ne i_2$ and $\gamma_{i_{1},1}\gamma_{i_{2},1}\ne 0$,
then we have $c(\gamma,0)=0$. 
\end{sublema}
\begin{pf}
We have 
\[
\deg (t^{\gamma})\ge (k+4)\frac{a_{l}-1}{a_{l}}>2, 
\]
where $a_{l}=\min \{a_{i_{1}},a_{i_{2}}\}$.
Hence we have $c(\gamma,0)=0$.   
\qed
\end{pf}

\vspace{5pt}
\noindent
\underline {\it Step 3. $a_{1}=a_{2}=2$ and $a_3\ge 3$, i.e., $A$ is a semi--general multiplet.}

\vspace{5pt}
If $t^{\gamma}$ has, as a factor, $t_{i_3,j_3}$ for $i_3\ge 3$ and $j_3\ge 1$, 
we have $c(\gamma,0)=0$ by previous arguments in Step 1 and Step 2. 
Therefore we only have to show the following Sublemma~\ref{sl6-l2-sep}.
\vspace{-5pt}
\begin{sublemd}[Step 3]\label{sl6-l2-sep}
Suppose that $A$ is a smi--general multiplet.
Then we have 
\[
\begin{cases}
{\rm (i)} \ \ c(2e_{1,1}+2e_{2,1},0)=0,\\
{\rm (ii)} \ \ c(3e_{1,1}+e_{2,1},0)=c(3e_{2,1}+e_{1,1},0)=0.
\end{cases}
\]
\end{sublemd}
\begin{pf}
Note that $c(\gamma,1)\neq 0$ with $|\gamma|=3$
if and only if $\gamma=\sum_{i=1}^{r}e_{i,1}$ by Proposition~\ref{lem3.1}, and that
$c(\gamma,0)=0$ if $\gamma \in \ZZ^{\mu_{A}-2}$ is a non--negative element such that 
$|\gamma|=4$ and $\gamma -e_{i_{1},j_{1}}-e_{i_3,j_{3}}\ge 0$ for $i_{1}=1,2$ and $i_3\ge 3$
by Lemma~\ref{wsep}.

First, we shall show that $c(2e_{1,1}+2e_{2,1},0)=0$. 
We shall calculate the coefficient of the term $t_{2,1}^{2}(\prod_{i=3}^{r}t_{i,1})e^{t_{\mu_A}}$ in $WDVV((1,1),(2,1),\mu_A,\mu_A)$.
Then we have
\[
4\cdot c(2e_{1,1}+2e_{2,1},0)\cdot 2\cdot 1=0. 
\]
Hence we have $c(2e_{1,1}+2e_{2,1},0)=0$.

Next, we shall show that $c(3e_{1,1}+e_{2,1},0)=c(3e_{2,1}+e_{1,1},0)=0$.
We shall calculate the coefficient of the term $t^{2}_{1,1}(\prod_{i=3}^{r}t_{i,1})e^{t_{\mu_A}}$ in $WDVV((1,1),(1,1),\mu_A,\mu_A)$. 
Then we have
\[
6\cdot c(3e_{1,1}+e_{2,1},0)\cdot a_{2}\cdot c(e_{1,1}+e_{2,1}+e_{3,1},1)=0. 
\]
Thus we have $c(3e_{1,1}+e_{2,1},0)=0$.
The same argument shows $c(3e_{2,1}+e_{1,1},0)=0$. 
\qed
\end{pf}
Therefore we have Theorem~\ref{sep}.
\qed
\end{pf}

For a non--general multiplet $A$, we have the following conjecture:
\begin{conj}
For each non--general multiplet $A$, there exsists a Frobenius structure which satisfies the conditions
{\rm (i), (ii), (iii), (v), (vi)} in Theorem~\ref{first}
and does not satisfy the condition {\rm (iv)}.
\end{conj}
As we remarked at the former part of this section, for $A=(2,2,2)$, one can obtain the Frobenius potential which 
satisfies the conditions in Theorem~\ref{first} except for the condition {\rm (iv)} by easy calculation.
\section{The Gromov-Witten Theory for Orbifold Projective Lines} 

Let $r\ge 3$ be a positive integer.
Let $A = (a_{1},\ldots,a_{r})$ be a multiplet of positive integers and 
$\Lambda=(\lambda_1,\ldots,\lambda_r)$ a multiplet of pairwise distinct elements of $\PP^{1}(\CC)$ normalized such that 
$\lambda_1=\infty$, $\lambda_2=0$ and $\lambda_3=1$. 
Following Geigle--Lenzing (cf. Section~1.1 in \cite{gl:1}), we shall introduce an orbifold projective line.
First, we prepare some notations.
\begin{defn}\label{orb proj line}
Let $r$, $A$ and $\Lambda$ be as above.
\begin{enumerate}
\item
Define a ring $R_{A,\Lambda}$ by 
\begin{subequations}
\begin{equation} 
R_{A,\Lambda}:=\CC[X_1,\dots,X_r]\left/I_{\Lambda}\right.,
\end{equation}
where  $I_\Lambda$ is an ideal generated by $r-2$ homogeneous polynomials
\begin{equation}
X_i^{a_i}-X_2^{a_2}+\lambda_i X_1^{a_1},\quad i=3,\dots, r.
\end{equation}
\end{subequations}
\item
Denote by $L_A$ an abelian group generated by $r$-letters $\vec{X_i}$, $i=1,\dots ,r$ defined as the quotient 
\begin{subequations}
\begin{equation}
L_A:=\bigoplus_{i=1}^r\ZZ\vec{X}_i\left/M_A\right. ,
\end{equation}
where  $M_A$ is the subgroup generated by the elements
\begin{equation}
a_i\vec{X}_i-a_j\vec{X}_j,\quad 1\le i<j\le r.
\end{equation}
\end{subequations}
\end{enumerate}
\end{defn}
We then consider the following quotient stack$:$
\begin{defn}\label{defn:gl}
Let $r$, $A$ and $\Lambda$ be as above.
Define a stack $\PP^1_{A,\Lambda}$ by
\begin{equation}
\PP^1_{A,\Lambda}:=\left[\left({\rm Spec}(R_{A,\Lambda})\backslash\{0\}\right)/{\rm Spec}({\CC L_A})\right],
\end{equation}
which is called the {\it orbifold projective line} of type $(A,\Lambda)$. 
\end{defn}
An orbifold projective line of type $(A,\Lambda)$ is a Deligne--Mumford stack whose coarse moduli space is 
a smooth projective line $\PP^1$.
The orbifold cohomology group of $\PP^{1}_{A,\Lambda}$ is, as a vector space, just the
singular cohomology group of the inertia orbifold: 
\[
\I\PP^{1}_{A,\Lambda}=\PP^{1}_{A,\Lambda} \bigsqcup \bigsqcup_{1\le i\le r} \left(\bigsqcup^{a_{i}-1}_{j=1} \left(B(\ZZ/a_i \ZZ)\right)_{j}\right)
 \]
where $\left(B(\ZZ/a_i \ZZ)\right)_{j}:=B(\ZZ/a_i \ZZ)$. The age associated to the component $\PP^{1}_{A,\Lambda}$
is $0$ and the age associated to $\left(B(\ZZ/a_i \ZZ)\right)_{j}$ is $j/a_{i}$.
The orbifold Poincar\'{e} pairing is given by twisting the usual Poincar\'{e} pairing:
\[
\displaystyle
\int_{\PP^{1}_{A,\Lambda}} \alpha \cup_{orb} \beta := \int_{\I\PP^{1}_{A,\Lambda}} \alpha \cup I \beta,
\]
where $I$ is the involution defined in \cite{agv:1, cr:1}.
Then we have the following:
\begin{lema}\label{gw1}
We can choose a $\QQ$-basis $1=\Delta_{1},\Delta_{1,1},\dots,\Delta_{i,j},\dots,\Delta_{r,a_r-1},\Delta_{\mu_A}$ of 
the orbifold cohomology group $H^{*}_{orb}(\PP^{1}_{A,\Lambda},\QQ)$
such that
\[
H^{0}_{orb}(\PP^{1}_{A,\Lambda},\QQ )\simeq \QQ \Delta_{1}, \ \  
\Delta_{i,j}\in H^{2\frac{j}{a_{i}}}_{orb}(\PP^{1}_{A,\Lambda},\QQ ), \ \
H^{2}_{orb}(\PP^{1}_{A,\Lambda},\QQ )\simeq \QQ \Delta_{\mu_{A}}
\]
and
\[
\displaystyle\int_{\PP^{1}_{A,\Lambda}} \Delta_{1} \cup_{orb} \Delta_{\mu_A} =1, \ \
\displaystyle\int_{\PP^{1}_{A,\Lambda}} \Delta_{i,j} \cup_{orb} \Delta_{k,l} =
\begin{cases}
&\frac{1}{a_{i}} \ \ \text{if} \ \ k=i, \ l=a_{i}-j \\
&0 \ \ \ \text{otherwise} .
\end{cases}
\]
\end{lema}
\begin{pf}
The decomposition of $H^{*}_{orb}(\PP^{1}_{A,\Lambda},\CC)$ follows from the decomposition of the inertia orbifold
$\I\PP^{1}_{A,\Lambda}$. The latter assertion immediately follows from the definition of the orbifold Poincar\'{e} pairing.
\qed
\end{pf}
Denote by $t_{1},t_{1,1},\dots ,t_{i,j},\dots , t_{r,a_r-1},t_{\mu_A}$ the dual coordinates of the $\QQ$-basis 
$\Delta_{1}$, $\Delta_{1,1},\dots,\Delta_{i,j},\dots,\Delta_{r,a_r-1},\Delta_{\mu_A}$ of  $H^{*}_{orb}(\PP^{1}_{A,\Lambda},\QQ)$ in Lemma~\ref{gw1}. 
Consider a formal manifold $M$ whose structure sheaf $\O_M$ and tangent sheaf $\T_M$ are given by 
\begin{equation}
\O_M:=\CC((e^{t_{\mu_A}}))[[t_1,t_{1,1},\dots ,t_{i,j},\dots , t_{r,a_r-1}]],\quad \T_M:=H^*_{orb}(\PP^1_{A,\Lambda},\CC)\otimes_\CC\O_M,
\end{equation}
where $\CC((e^{t_{\mu_A}}))$ denotes the $\CC$-algebra of formal Laurent series in $e^{t_{\mu_A}}$.
The Gromov--Witten theory for orbifolds developed by Abramovich--Graber--Vistoli \cite{agv:1} and Chen--Ruan \cite{cr:1} 
gives us the following proposition. 
Note here that, by using the divisor axiom, it turns out that third derivatives of the genus zero Gromov--Witten potential $\F^{\PP^{1}_{A,\Lambda}}_{0}$ are elements of $\CC[[t_{1,1},\dots ,t_{i,j},\dots , t_{r,a_r-1}, q^{[\PP^1]}e^{t_{\mu_A}}]]$ and hence
they can be considered as elements of $\O_M$ by formally setting $q^{[\PP^1]}=1$.
\begin{prop}[\cite{agv:1, cr:1}]
There exists a structure of a formal Frobenius manifold of rank $\mu_A$ and dimension one 
on $M$ whose non--degenerate symmetric $\O_M$--bilinear form $\eta$ on $\T_M$ is given by the orbifold Poincar\'{e} pairing.
\end{prop}
\begin{pf}
See Theorem 6.2.1 of \cite{agv:1} and Theorem 3.4.3 of \cite{cr:1}.
\qed
\end{pf}

The following theorem is the main result in this section:
\begin{thm}\label{second}
The conditions in Theorem~\ref{first} are satisfied by the Frobenius structure constructed
from the Gromov--Witten theory for $\PP^1_{A,\Lambda}$. 
\end{thm}
We shall check the conditions in Theorem~\ref{first} one by one.
\subsection{Condition (i)}
It follows from Lemma~\ref{gw1}
that the unit vector field $e\in\T_M$ and the Euler vector field $E\in\T_M$ are given as
\[
e=\frac{\p}{\p t_1}, \ \ E=t_1\frac{\p}{\p t_1}+\sum_{i=1}^{r}\sum_{j=1}^{a_i-1}\frac{a_i-j}{a_i}t_{i,j}\frac{\p}{\p t_{i,j}}
+\chi_{A}\frac{\p}{\p t_{\mu_A}},
\]
which is the condition {\rm (i)}.
\subsection{Condition (ii)}
It is obvious from Lemma~\ref{gw1}.
\subsection{Condition (iii)}
The condition {\rm (iii)} follows from the divisor axiom and 
the definition of the genus zero potential $\F^{\PP^{1}_{A,\Lambda}}_{0}$.

\subsection{Condition (iv)}\label{con(4)}
The condition {\rm (iv)} is satisfied since the image of degree zero orbifold map with marked points on orbifold points 
on the source must be one of orbifold points on the target $\PP^1_{A,\Lambda}$.

\subsection{Condition (v)}

The orbifold cup product is the specialization of the quantum product
at $t_{1}=t_{1,1}=\dots=t_{r,a_r-1}=e^{t_{\mu_A}}=0$. 
Therefore, it turns out that the orbifold cup product can be determined by the degree zero
three point Gromov-Witten invariants.
\begin{lema}\label{gw5}
There is a $\CC$-algebra isomorphism between the orbifold cohomology ring $H^{*}_{orb}(\PP^{1}_{A,\Lambda},\CC)$ and
$\CC[x_1,x_2,\dots, x_r]\left/\left(x_ix_j, \ a_ix_i^{a_i}-a_jx_j^{a_j}
\right)_{1\le i\ne j\le r}\right.,$
where $\p/\p t_{i,j}$ are mapped to
$x^{j}_i$ for $i=1,\dots,r, j=1,\dots, a_{i}-1$ and $\p/\p t_{\mu_A}$ are mapped to $a_{1}x_{1}^{a_1}$.
\end{lema}
\begin{pf}
Under the same notation in Lemma~\ref{gw1},
the orbifold cup product is given as follows:
\[
\Delta_{\alpha}\cup_{orb}\Delta_{\beta}=\sum_{\delta}\left<\Delta_{\alpha},\Delta_{\beta},\Delta_{\gamma}\right>^{\PP^{1}_{A,\Lambda}}_{0,3,0}
\eta^{\gamma\delta} \Delta_{\delta},
\]
where we set $\eta^{\gamma\delta}$ as follows:
\[
\displaystyle(\eta^{\gamma\delta})=(\int_{\PP^{1}_{A,\Lambda}}\Delta_{\gamma}\cup_{orb}\Delta_{\delta})^{-1}.
\]
By the previous argument in Subsection \ref{con(4)}, 
we have
\[
\Delta_{i_1,j_1}\cup_{orb} \Delta_{i_2,j_2}=0 \ \text{if} \ i_1\ne i_2.
\] 
By the formula
\begin{align*}
 \int_{\PP^{1}_{A,\Lambda}} \Delta_{i_1,j_1} \cup_{orb} \Delta_{i_1,j'_1} \cup_{orb} \Delta_{i_1,j''_1}
=&\  \frac{1}{|\ZZ/a_{i_1}\ZZ|} \int_{pt} ev^{*}_1 (\Delta_{i_1,j_1})\cup ev^{*}_2 (\Delta_{i_1,j'_1})\cup ev^{*}_3 (\Delta_{i_1,j''_1})\\
=&\ 
\begin{cases}
\frac{1}{a_{i_1}} \ \ \ \text{if} \ j_1 +j'_1 +j''_1=a_{i_1}, \\
0\ \ \ \text{otherwise},
\end{cases}
\end{align*}
we have
\[
\Delta_{i_{1},j_{1}}\cup_{orb} \Delta_{i_1,j'_1}=\Delta_{i_1,j_1+j'_1} \ \text{if} \ j_1+j'_1\le a_{i_1}-1,
\]
and hence
\[
\displaystyle\Delta_{i_1,1}^{p}:=\underbrace{\Delta_{i_1,1}\cup_{orb}\cdots \cup_{orb} \Delta_{i_1,1}}_{a_{i_1} \ times}=
\frac{1}{a_{i_1}}\Delta_{\mu_A}.
\]
Therefore we have Lemma~\ref{gw5}.
\qed
\end{pf}

\begin{lema}\label{gw6}
The term 
\[
\displaystyle\left(\prod_{i=1}^{r}t_{i,1}\right)e^{t_{\mu_A}}
\]
occurs with the coefficient $1$ in the $\F^{\PP^{1}_{A,\Lambda}}_{0}$.
\end{lema}
\begin{pf}
This lemma follows from the fact that the Gromov-Witten invariant counts the number of
orbifold maps from $\PP^{1}_{A,\Lambda}$ to $\PP^{1}_{A,\Lambda}$ of degree $1$ fixing $r$ marked $($orbifold$)$ points, 
which is exactly the identity map. 
\qed
\end{pf}
\section{Vanishing of higher degree correlators}

As we see in Subsection~\ref{reconst-s2}, the coefficients $c(e_{i,1}+e_{i,a_{i}-1}+e_{1},0)$, $c(e_{i,1}+e_{i,j-1}+e_{i,a_{i}-j},0)$, 
$c(\sum_{i=1}^{r}e_{i,1},1)$ or $c(2e_{i,1}+2e_{i,a_{i}-1},0)$ play important roles to reconstruct the Frobenius potential.
Indeed, in the works of Krawitz--Shen (\cite{ks:1}) and Li--Li--Saito--Shen (\cite{llss:1}), the similar type reconstuction theorems were proved and
the authors also used them as the initial data in order to reconstruct the Frobenius potential. 
The coefficient $c(\sum_{i=1}^{r}e_{i,1},1)$ corresponds to a certain degree one $r$-points correlator in the orbifold Gromov--Witten theory of $\PP^{1}_{A,\Lambda}$
and hence can be determined easily by considering its geometric meaning. However, in other theory like the 
invariant theory of extended cuspidal Weyl groups which is expected as a mirror partner of $\PP^{1}_{A,\Lambda}$ (\cite{Sh-T:1}), the representation theoretic
meaning of this coefficient is not known at all and it is hard even to verify whether this coefficient is non--zero or not.  
On the other hand, it seems we are able to check other conditions in Theorem \ref{first} one way or another.
Therefore it is very interesting and important to investigate what happens on
the Frobenius potential which satisfies the conditions in Theorem \ref{first} except for the condition {\rm (vi)}.

For the Frobeius potential above, we have the following two cases analysis:
\begin{equation*}
c(\sum_{i=1}^{r}e_{i,1},1)=
\begin{cases}
{\rm (i)} \quad a\ne 0,\\
{\rm (ii)} \quad 0.
\end{cases}
\end{equation*}
For the case {\rm (i)}, we can obtain the following Frobenius potential $\F'_{A}$:
\begin{equation}
\F'_{A}=\sum_{\alpha} c(\alpha ,0) t^{\alpha}+\sum_{\alpha} a^{m}c(\alpha ,m) t^{\alpha}e^{mt_{\mu_A}},
\end{equation}
where coefficients $c(\alpha ,0)$ and $c(\alpha ,m)$ are same ones of the Frobenius potential in Theorem~\ref{first}.
This result can be shown easily by checking the each procedure in the proof of Theorem~\ref{first} carefully.
In particular, the Frobenius potential $\F'_A$ can be transformed to the one in Theorem~\ref{first} by the (invertible) coordinate change
$t_{\mu_A}\mapsto t_{\mu_A}-{\rm log}(a)$.

The second case {\rm (ii)} is much more serious since the condition {\rm (vi)} is heavily used for determining all other coefficients except
for $c(\alpha,0)$ with $|\alpha|=3$. Actually, a Frobenius potential under the case {\rm (ii)} is no longer determined uniquely even if $A$ is of type $ADE$.
For this problem, we have to consider another reasonable initial amount, i.e., 
the coefficient $c(2e_{i,1}+2e_{i,a_{i}-1},0)$ which corresponds to a degree zero $4$-points correlator in the orbifold Gromov--Witten theory.
By this reason, we shall consider the following two cases:
\begin{equation}
c(2e_{i,1}+2e_{i,a_{i}-1},0)=
\begin{cases}
{\rm (i)} \ \ -1/4a_{i}^{2} \ \ \text{if} \ \ a_{i}\ge 3, \ \ -1/96 \ \ \text{if} \ \ a_{i}=2, \\
{\rm (ii)} \ \ \text{otherwise}.
\end{cases}
\end{equation}

For the first case {\rm (i)}, we have the following proposition:
\begin{prop}\label{prop:coeff 0}
Assume that a Frobenius manifold $M$ of rank $\mu_A$ and dimension one with flat coordinates 
$(t_1,t_{1,1},\dots ,t_{i,j},\dots ,t_{r,a_r-1},t_{\mu_A})$ satisfies the following conditions$:$
\begin{enumerate}
\item 
The unit vector field $e$ and the Euler vector field $E$ are given by
\[
e=\frac{\p}{\p t_1},\ E=t_1\frac{\p}{\p t_1}+\sum_{i=1}^{r}\sum_{j=1}^{a_i-1}\frac{a_i-j}{a_i}t_{i,j}\frac{\p}{\p t_{i,j}}
+\chi_A\frac{\p}{\p t_{\mu_A}}.
\]
\item 
The non--degenerate symmetric bilinear form $\eta$ on $\T_M$ satisfies
\begin{align*}
&\ \eta\left(\frac{\p}{\p t_1}, \frac{\p}{\p t_{\mu_A}}\right)=
\eta\left(\frac{\p}{\p t_{\mu_A}}, \frac{\p}{\p t_1}\right)=1,\\ 
&\ \eta\left(\frac{\p}{\p t_{i_1,j_1}}, \frac{\p}{\p t_{i_2,j_2}}\right)=
\begin{cases}
\frac{1}{a_{i_1}}\quad i_1=i_2\text{ and }j_2=a_{i_1}-j_1,\\
0 \quad \text{otherwise}.
\end{cases}
\end{align*}
\item 
The Frobenius potential $\F$ satisfies $E\F|_{t_{1}=0}=2\F|_{t_{1}=0}$,
\[
\left.\F\right|_{t_1=0}\in\CC\left[[t_{1,1}, \dots, t_{1,a_1-1}, 
\dots, t_{i,j},\dots, t_{r,1}, \dots, t_{r,a_r-1},e^{t_{\mu_A}}]\right].
\]
\item Assume the condition {\rm (iii)}. we have
\begin{equation*}
\F|_{t_1=e^{t_{\mu_A}}=0}=\sum_{i=1}^{r}\G^{(i)}, \quad \G^{(i)}\in \CC[[t_{i,1},\dots, t_{i,a_i-1}]],\ i=1,\dots,r.
\end{equation*}
\item 
Assume the condition {\rm (iii)}. In the frame $\frac{\p}{\p t_1}, \frac{\p}{\p t_{1,1}},\dots, 
\frac{\p}{\p t_{r,a_r-1}},\frac{\p}{\p t_{\mu_A}}$ of $\T_M$,
the product $\circ$ can be extended to the limit $t_1=t_{1,1}=\dots=t_{r,a_r-1}=e^{t_{\mu_A}}=0$.
The $\CC$-algebra obtained in this limit is isomorphic to
\[
\CC[x_1,x_2,\dots, x_r]\left/\left(x_ix_j, \ a_ix_i^{a_i}-a_jx_j^{a_j}
\right)_{1\le i\ne j\le r}\right.,
\]
where $\p/\p t_{i,j}$ are mapped to
$x^{j}_i$ for $i=1,\dots,r, j=1,\dots, a_{i}-1$ and $\p/\p t_{\mu_A}$ are mapped to $a_{1}x_{1}^{a_1}$.
\item The term 
\[
\displaystyle\left(\prod_{i=1}^{r}t_{i,1}\right)e^{t_{\mu_A}}
\]
occurs with the coefficient $0$ in $\F$. 
\item 
The term $t_{i,1}^{2}t_{i,a_{i}-1}^{2}$ in $\F$ occurs with the coefficient 
\[
\begin{cases}
\displaystyle
-1/96 \quad \text{if} \quad a_{i}=2,\\
\displaystyle
-1/4a_{i}^{2} \quad \text{if} \quad a_{i}\ge3.
\end{cases}
\]
\end{enumerate}
Then any term $t^{\alpha}e^{mt_{\mu_A}}$ for $m\ge 1$ occurs with the coefficient $0$ in $\F$.
\end{prop}

\begin{pf}
By the same argument in Lemma~\ref{lem1-lem3.1} and Lemma~\ref{lem2-lem3.1},
we have $c(\gamma,1)=0$ for $|\gamma|\le r$.
We shall show Proposition~\ref{prop:coeff 0} by the induction on the total oreder defined in Subsection~\ref{reconst-s2}.
\begin{lem}\label{lem0-mreconst:0}
Assume that $c(\alpha, n)=0$ for $(|\alpha|,n)\prec (0,m)$.
Then we have $c(0,m)=0$
\end{lem}

\begin{pf}
We shall calculate the coefficient of the term  
$e^{mt_{\mu_A}}$
in $WDVV((i,1),(i,a_{i}-1),\mu_A,\mu_A)$. Then we have
\[
c(e_{i,1}+e_{i,a_{i}-1}+e_{1},0)\cdot 1 \cdot m^{3}\cdot c(0,m)=0
\]
since $c(\alpha, n)=0$ for $(|\alpha|,n)\prec (0,m)$. Therefore we have $c(0,m)=0$.
\qed
\end{pf}
Next, we shall split the second step of the induction into following three cases. 
\begin{lem}[Case 1]\label{lem1-mreconst:0}
Assume that $c(\alpha, n)=0$ for $(|\alpha|,n)\prec (k+1,m)$.
If a non--negative element $\gamma \in \ZZ^{\mu_{A}-2}$ satisfies that $|\gamma|=k+1$
and $\gamma-e_{i,j}\ge 0$ for some $j$ such that $j\ge 2$,
then we have $c(\gamma,m)=0$.
\end{lem}
\begin{pf}
We shall calculate the coefficient of the term  
$t^{\gamma-e_{i,j}}e^{mt_{\mu_A}}$
in the WDVV equation $WDVV((i,1),(i,j-1),\mu_A,\mu_A)$. Then we have
\[
s_{1,j-1,a_{i}-j}\cdot c(e_{i,1}+e_{i,j-1}+e_{i,a_{i}-j},0)\cdot a_{i} \cdot m^{2}\cdot \gamma_{i,j} \cdot c(\gamma ,m)=0
\]
since we have $c(\alpha, n)=0$ for $(|\alpha|,n)\prec (k+1,m)$. Hence we have $c(\gamma, m)=0$.
\qed
\end{pf}
\begin{lem}[Case 2]\label{lem2-mreconst:0}
Assume that $c(\alpha, n)=0$ for $(|\alpha|,n)\prec (k+1,m)$.
If a non--negative element $\gamma \in \ZZ^{\mu_{A}-2}$ satisfies that $|\gamma|=k+1$
and $\gamma =\sum_{k=1}^{r}\gamma_{k,1}e_{k,1}$ for some  $\gamma_{1,1},\dots,\gamma_{r,1}$ such that $\prod_{k=1}^{r}\gamma_{k,1}\ne 0$,
then we have $c(\gamma ,m)=0$.
\end{lem}
\begin{pf}
We shall calculate the coefficient of the term $(\prod_{k=1}^{r}t_{k,1}^{\gamma_{k,1}})e^{mt_{\mu_A}}$
in the WDVV equation $WDVV((i,1),(i,a_{i}-1),\mu_A,\mu_A)$. Then we have, if $a_{i}\ge 3$,
\begin{equation*}
\{c(e_{i,1}+e_{i,a_{i}-1}+e_{1},0)\cdot m^{3} +4\cdot c(2e_{i,1}+2e_{i,a_{i}-1},0)\cdot a_{i} 
\cdot m^{2}\cdot \gamma_{i,1}\} \cdot c(\gamma ,m)=0 
\end{equation*}
and, if $a_{i}=2$,
\begin{equation*}
\{ 2c(2e_{i,1}+e_{1},0)\cdot m^{3} +24\cdot c(4e_{i,1},0)\cdot 2\cdot m^{2} 
\cdot \gamma_{i,1}\} \cdot c(\gamma ,m)=0 
\end{equation*}
since we have $c(\alpha, n)=0$ for $(\alpha,n)\prec (k+1,m)$. If $\gamma_{i,1}\neq m$ for some $i$, this lemma holds.
If $\gamma_{i,1}=m$ for all $i$, the degree $\deg ((\prod_{k=1}^{r}t_{k,1}^{\gamma_{k,1}})e^{mt_{\mu_A}})$ 
is greater than $2$ except for the case $m=1$ and hence this lemma also holds for this case. 
\qed
\end{pf}
\begin{lem}[Case 3]\label{lem3-mreconst:0}
Assume that $c(\alpha, n)=0$ for $(|\alpha|,n)\prec (k+1,m)$.
If a non--negative element $\gamma \in \ZZ^{\mu_{A}-2}$ satisfies that $|\gamma|=k+1$ and 
$\gamma =\sum_{k=1}^{r}\gamma_{k,1}e_{k,1}$ for some  $\gamma_{1,1},\dots,\gamma_{r,1}$ such that $\prod_{k=1}^{r}\gamma_{k,1}= 0$,
then we have $c(\gamma ,m)=0$ with $(|\alpha|,n)\prec (k+1,m)$.
\end{lem}
\begin{pf}
Assume that $\gamma_{i,1}=0$.
We shall calculate the coefficient of the term $(\prod_{k=1}^{r}t_{k,1}^{\gamma_{k,1}})e^{mt_{\mu_A}}$
in the equation $WDVV((i,1),(i,a_{i}-1),\mu_A,\mu_A)$. Then we have
\[
c(e_{1}+e_{i,a_{i}-1}+e_{i,1},0)\cdot m^{3} \cdot c(\gamma ,m)=0
\]
since $c(\alpha, n)=0$ for $(|\alpha|,n)\prec (k+1,m)$. Then we have $c(\gamma, m)=0$
\qed
\end{pf}
We have finished the proof of this proposition.
\qed
\end{pf}
\begin{rem}
Proposition~\ref{prop:coeff 0} would have an application to show the uniqueness of the Frobenius manifold $M_{\widehat{W}_A}$ constructed from the invariant theory of
an extended cuspidal Weyl group in a further joint work \cite{Sh-T:1} and an isomorphism of Frobenius manifolds between $M_{\widehat{W}_A}$ and 
$M^{GW}_{\PP^{1}_{A,\Lambda}}$.  
\end{rem}

For the case {\rm (ii)}, we have the following proposition: 
\begin{prop}\label{prop:vanish-higher}
Assume that $\chi_A\ge 0$ and a Frobenius manifold $M$ of rank $\mu_A$ and dimension one with flat coordinates 
$(t_1,t_{1,1},\dots ,t_{i,j},\dots ,t_{r,a_r-1},t_{\mu_A})$ satisfies the following conditions$:$
\begin{enumerate}
\item 
The unit vector field $e$ and the Euler vector field $E$ are given by
\[
e=\frac{\p}{\p t_1},\ E=t_1\frac{\p}{\p t_1}+\sum_{i=1}^{r}\sum_{j=1}^{a_i-1}\frac{a_i-j}{a_i}t_{i,j}\frac{\p}{\p t_{i,j}}
+\chi_A\frac{\p}{\p t_{\mu_A}}.
\]
\item 
The non--degenerate symmetric bilinear form $\eta$ on $\T_M$ satisfies
\begin{align*}
&\ \eta\left(\frac{\p}{\p t_1}, \frac{\p}{\p t_{\mu_A}}\right)=
\eta\left(\frac{\p}{\p t_{\mu_A}}, \frac{\p}{\p t_1}\right)=1,\\ 
&\ \eta\left(\frac{\p}{\p t_{i_1,j_1}}, \frac{\p}{\p t_{i_2,j_2}}\right)=
\begin{cases}
\frac{1}{a_{i_1}}\quad i_1=i_2\text{ and }j_2=a_{i_1}-j_1,\\
0 \quad \text{otherwise}.
\end{cases}
\end{align*}
\item 
The Frobenius potential $\F$ satisfies $E\F|_{t_{1}=0}=2\F|_{t_{1}=0}$,
\[
\left.\F\right|_{t_1=0}\in\CC\left[[t_{1,1}, \dots, t_{1,a_1-1}, 
\dots, t_{i,j},\dots, t_{r,1}, \dots, t_{r,a_r-1},e^{t_{\mu_A}}]\right].
\]
\item Assume the condition {\rm (iii)}. we have
\begin{equation*}
\F|_{t_1=e^{t_{\mu_A}}=0}=\sum_{i=1}^{r}\G^{(i)}, \quad \G^{(i)}\in \CC[[t_{i,1},\dots, t_{i,a_i-1}]],\ i=1,\dots,r.
\end{equation*}
\item 
Assume the condition {\rm (iii)}. In the frame $\frac{\p}{\p t_1}, \frac{\p}{\p t_{1,1}},\dots, 
\frac{\p}{\p t_{r,a_r-1}},\frac{\p}{\p t_{\mu_A}}$ of $\T_M$,
the product $\circ$ can be extended to the limit $t_1=t_{1,1}=\dots=t_{r,a_r-1}=e^{t_{\mu_A}}=0$.
The $\CC$-algebra obtained in this limit is isomorphic to
\[
\CC[x_1,x_2,\dots, x_r]\left/\left(x_ix_j, \ a_ix_i^{a_i}-a_jx_j^{a_j}
\right)_{1\le i\ne j\le r}\right.,
\]
where $\p/\p t_{i,j}$ are mapped to
$x^{j}_i$ for $i=1,\dots,r, j=1,\dots, a_{i}-1$ and $\p/\p t_{\mu_A}$ are mapped to $a_{1}x_{1}^{a_1}$.
\item The term 
\[
\displaystyle\left(\prod_{i=1}^{r}t_{i,1}\right)e^{t_{\mu_A}}
\]
occurs with the coefficient $0$ in $\F$. 
\end{enumerate}
Then any term $t^{\alpha}e^{mt_{\mu_A}}$ for $m\ge 1$ occurs with the coefficient $0$ in $\F$.
\end{prop}
\begin{pf}
Put $\gamma :=\sum_{k=1}^{r}\gamma_{k,1}e_{k,1}$ for some  $\gamma_{1,1},\dots,\gamma_{r,1}$ such that $\prod_{k=1}^{r}\gamma_{k,1}\ne 0$.
If $m\ge 2$, then we have $c(\gamma, m)=0$ since we have $\chi_A\ge 0$ and the following inequation:
\begin{equation}
{\deg}(t^{\gamma}e^{mt_{\mu_A}})=\sum^{r}_{i=1}\gamma_{i,1}\frac{a_{i}-1}{a_{i}}+m\chi_{A} >\sum^{r}_{i=1}\frac{a_{i}-1}{a_{i}}+\chi_{A}=2.
\end{equation}
Same arguments in Lemma~\ref{lem0-mreconst:0}, Lemma~\ref{lem1-mreconst:0} and Lemma~\ref{lem3-mreconst:0} shows this proposition.
\qed
\end{pf}
Same statement as in Proposition~\ref{prop:vanish-higher} would hold for $\chi_A<0$ by the exactly the same arguments 
if there were not positive integers $m$ and $\gamma_{i,1}$ ($i=1,\dots,r$) satisfying
the following equations:
\begin{equation}
c(e_{i,1}+e_{i,a_{i}-1}+e_{1},0)\cdot m^{3} +4\cdot c(2e_{i,1}+2e_{i,a_{i}-1},0)\cdot a_{i} 
\cdot m^{2}\cdot \gamma_{i,1}=0 \quad \text{if} \ \ a_{i}\ge 3,
\end{equation}
\begin{equation}
2c(2e_{i,1}+e_{1},0)\cdot m^{3} +24\cdot c(4e_{i,1},0)\cdot 2\cdot m^{2} 
\cdot \gamma_{i,1}=0 \quad \text{if} \ \ a_{i}=2,
\end{equation}
\begin{equation}
\sum^{r}_{i=1}\gamma_{i,1}\frac{a_{i}-1}{a_{i}}+m\chi_{A}=2.
\end{equation}
However, for the cases $\chi_A<0$, one can choose $4$-points correlators $c(2e_{i,1}+2e_{i,a_{i}-1},0)$ such that
there exist positive integers $m$ and $\gamma_{i,1}$ ($i=1,\dots,r$) satisfying the above equations.
By the computer experiment, we conjecture the following: 
\begin{conj}\label{conj}
There exists a Frobenius manifold $M$ of rank $\mu_A$ and dimension one with flat coordinates 
$(t_1,t_{1,1},\dots ,t_{i,j},\dots ,t_{r,a_r-1},t_{\mu_A})$ satisfies the following conditions$:$
\begin{enumerate}
\item 
The unit vector field $e$ and the Euler vector field $E$ are given by
\[
e=\frac{\p}{\p t_1},\ E=t_1\frac{\p}{\p t_1}+\sum_{i=1}^{r}\sum_{j=1}^{a_i-1}\frac{a_i-j}{a_i}t_{i,j}\frac{\p}{\p t_{i,j}}
+\chi_A\frac{\p}{\p t_{\mu_A}}.
\]
\item 
The non--degenerate symmetric bilinear form $\eta$ on $\T_M$ satisfies
\begin{align*}
&\ \eta\left(\frac{\p}{\p t_1}, \frac{\p}{\p t_{\mu_A}}\right)=
\eta\left(\frac{\p}{\p t_{\mu_A}}, \frac{\p}{\p t_1}\right)=1,\\ 
&\ \eta\left(\frac{\p}{\p t_{i_1,j_1}}, \frac{\p}{\p t_{i_2,j_2}}\right)=
\begin{cases}
\frac{1}{a_{i_1}}\quad i_1=i_2\text{ and }j_2=a_{i_1}-j_1,\\
0 \quad \text{otherwise}.
\end{cases}
\end{align*}
\item 
The Frobenius potential $\F$ satisfies $E\F|_{t_{1}=0}=2\F|_{t_{1}=0}$,
\[
\left.\F\right|_{t_1=0}\in\CC\left[[t_{1,1}, \dots, t_{1,a_1-1}, 
\dots, t_{i,j},\dots, t_{r,1}, \dots, t_{r,a_r-1},e^{t_{\mu_A}}]\right].
\]
\item Assume the condition {\rm (iii)}. we have
\begin{equation*}
\F|_{t_1=e^{t_{\mu_A}}=0}=\sum_{i=1}^{r}\G^{(i)}, \quad \G^{(i)}\in \CC[[t_{i,1},\dots, t_{i,a_i-1}]],\ i=1,\dots,r.
\end{equation*}
\item 
Assume the condition {\rm (iii)}. In the frame $\frac{\p}{\p t_1}, \frac{\p}{\p t_{1,1}},\dots, 
\frac{\p}{\p t_{r,a_r-1}},\frac{\p}{\p t_{\mu_A}}$ of $\T_M$,
the product $\circ$ can be extended to the limit $t_1=t_{1,1}=\dots=t_{r,a_r-1}=e^{t_{\mu_A}}=0$.
The $\CC$-algebra obtained in this limit is isomorphic to
\[
\CC[x_1,x_2,\dots, x_r]\left/\left(x_ix_j, \ a_ix_i^{a_i}-a_jx_j^{a_j}
\right)_{1\le i\ne j\le r}\right.,
\]
where $\p/\p t_{i,j}$ are mapped to
$x^{j}_i$ for $i=1,\dots,r, j=1,\dots, a_{i}-1$ and $\p/\p t_{\mu_A}$ are mapped to $a_{1}x_{1}^{a_1}$.
\item The term 
\[
\displaystyle\left(\prod_{i=1}^{r}t_{i,1}\right)e^{t_{\mu_A}}
\]
occurs with the coefficient $0$ in $\F$. 
\item A term
\[
\displaystyle t^{\gamma}e^{mt_{\mu_A}}
\]
occurs a non--zero coefficient in $\F$ for some non--negative $\gamma$ and $m\ge 2$.
\end{enumerate}
\end{conj}

\end{document}